\documentclass[11pt, letterpaper]{amsart}

\usepackage{amsmath, amssymb, amsthm, amsfonts}
\usepackage{color}
\usepackage{hyperref}
\hypersetup{colorlinks=true,allcolors=blue}

\newcommand{\R}{\mathbb{R}}
\newcommand{\N}{\mathbb{N}}

\newcommand{\uhat}{\widehat{u}}
\newcommand{\F}{\mathcal{F}}
\newcommand{\half}{\frac{1}{2}}
\newcommand{\norm}[1]{\|#1\|}
\newcommand{\abs}[1]{\left\lvert#1\right\rvert}
\newcommand{\jap}[1]{\left\langle#1\right\rangle}

\newtheorem{thm}{Theorem}[section]
\newtheorem{lemma}[thm]{Lemma}
\newtheorem{cor}[thm]{Corollary}

\newtheorem{propn}[thm]{Proposition}
\newtheorem{rem}[thm]{Remark}

\author[D. Bhattacharya]{Debdeep Bhattacharya}
\address{Department of Mathematics\\George Washington University, Washington, DC, USA}
\curraddr{}
\email{debdeepbh@gwu.edu}
\thanks{} 

\author[L. G. Farah]{Luiz Gustavo Farah}
\address{Department of Mathematics\\Universidade Federal de Minas Gerais\\Belo Horizonte, Brazil}
\curraddr{}
\email{farah@mat.ufmg.br}
\thanks{}

\author[S. Roudenko]{Svetlana Roudenko}
\address{Department of Mathematics \& Statistics\\Florida International University,  Miami, FL, USA}
\curraddr{}
\email{sroudenko@fiu.edu}
\thanks{}



\title[Global well-posedness in 2d modified ZK]{Global well-posedness for low regularity data\\ in the 2d modified Zakharov-Kuznetsov equation}

\begin{document}

\begin{abstract}
We consider the modified Zakharov-Kuznetsov (mZK) equation in two
space dimensions in both focusing and defocusing cases. Using the $I$-method, we prove the global well-posedness of the $H^s$ solutions for $s>\frac{3}{4}$ for any data in the defocusing case and under the assumption that the mass of the initial data is less than the mass of the ground state solution of $\Delta \varphi - \varphi + \varphi^3 = 0$ in the focusing case. This improves the global well-posedness result of Linares and Pastor \cite{LP}. 
\end{abstract}

\maketitle 

\section{Introduction}

We consider the two-dimensional initial value problem (IVP)
\begin{equation}
\begin{cases}
v_t + \partial_x (\Delta v)  + \sigma \partial_x (v^{3}) = 0 \text{, ~ ~ ~} (x,y) \in \R^2, t >0\\
v(x,y,0)= v_0(x,y),
\end{cases}
\label{agzk}
\end{equation}
where $v$ is a real-valued function, $\Delta = \partial_x^2 + \partial_y^2$ is the Laplacian operator in 2d and $\sigma = \pm1$ denotes the focusing (plus sign) and defocusing (minus sign) cases of the equation. 
This equation is a modification (thus, the name) of the standard Zakharov-Kuznetsov (ZK) equation 
$$
v_t + \partial_x(\Delta v + v^2) = 0,
$$
introduced in 3d by Zakharov and Kuznetsov \cite{ZK} to model the propagation of nonlinear ion-acoustic waves in magnetized plasma.

The mZK equation has been extensively studied in recent years, we recall relevant well-posedness results  in 2d.  
Biagioni and Linares \cite{BP} studied the local well-posedness of solutions to mZK with data in $H^1(\R^2)$.
Linares and Pastor \cite{LPlwp} proved the local well-posedness in $H^s(\R^2)$ for $s>\frac{3}{4}$ and they also showed the ill-posedness (non-uniform data to solution map) for $s \leq 0$, so one can not expect well-posedness in the critical space $L^2(\R^2)$.  
Ribaud and Vento \cite{RV} proved local well-posedness in $H^s(\R^2)$
for $s>\frac{1}{4}$, which is currently the best result known for the
local well-posedness for the 2d mZK equation. 
Linares and Pastor \cite{LP} proved the global well-posedness in
$H^s(\R^2)$ for $s> \frac{53}{63}$ with an additional assumption on the size of the initial data (related to the ground state) in the focusing case.  
In this paper, we use the I-method and obtain the global well-posedness in $H^s(\R^2)$ space for $s >\frac{3}{4}$, thus, improving the result of Linares and Pastor \cite{LP}.  

During their lifespan, solutions $u(x,y,t)$ to the mZK equation \eqref{agzk} conserve the mass $\mathcal{M}$ and energy $\mathcal{E}$: 
\begin{equation}
  \mathcal{M}[v(t)] = \int_{\R^2}^{}v^2(x,y,t) dxdy \equiv \mathcal{M}[v_0]
  \label{amass}
\end{equation}
and
\begin{equation}
  \mathcal{E}[v(t)] = \frac{1}{2}\int_{\R^2} |\nabla v(x,y,t)|^2
  dxdy - \frac{\sigma }{4} \int_{\R^2} v^{4}(x,y,t) dxdy \equiv \mathcal{E}[v_0].
  \label{aenergy}
\end{equation}

Note that if $v$ is a solution to \eqref{gzk}, then so is $v_\lambda$, its rescaled version
\begin{equation*}
v_\lambda(x,y,t) := \lambda v\left(\lambda x, \lambda y, \lambda^3 t\right).
\end{equation*}
The equation \eqref{agzk} is referred to as the $L^2$-critical (or mass-critical), since its $L^2$ norm is invariant under this rescaling. As with most focusing $L^2$-critical equations, solutions may blow up in finite time, which is also the case for this equation, see \cite{FHRY}, and thus, in order to consider the global well-posedness we have to put a restriction on the size of the initial data. For that we recall the notion of the ground state or traveling waves. 
Let $\varphi$ be the only radial and positive solution of
\begin{align} \label{GSeq}
\Delta \varphi -\varphi + \varphi^{3} = 0, \quad (x,y)\in \mathbb{R}^2.
\end{align}
Then
$$
v(x,y,t)=\varphi_c(x-ct,y)
$$
is a solution of the focusing mZK equation which travels only in the $x$-direction.  Here, $\varphi_c$ is the dilation of $\varphi$ given by
$$
\varphi_c(x,y)=\sqrt{c} \, \varphi(\sqrt{c}x,\sqrt{c}y)
$$
and solves the equation $\Delta \varphi_c -c\varphi_c + \varphi_c^{3} = 0$.

The existence of solutions of the equation \eqref{GSeq} in 2d was considered by Berestycki, Gallou\"et and  Kavian \cite{BGK83}, see  Strauss \cite{Sr77}, Berestycki and Lions \cite{BLi83}, Berestycki, Lions and Peletier \cite{BLP81} for the existence in other dimensions.  Regarding the uniqueness, Kwong \cite{Kwong} showed that radial and positive solutions are unique.

The function $\varphi$, also known as the \textit{ground state}, is related to the following sharp Gagliardo-Nirenberg inequality (see Weinstein \cite{Weinstein})
\begin{equation}\label{GNcrit}
\| {f} \|^4_{L^4(\R^2)} \le \frac{2}{\quad \| {\varphi} \|^2_{L^2(\R^2)}}\, \|{f} \|^2_{L^2(\R^2)} \, \| {\nabla f}\|^2_{L^2(\R^2)}.
\end{equation}

Recall that from the definition of the energy \eqref{aenergy} in the defocusing case $\sigma=-1$, we immediately have $\| {\nabla v(t)}\|^2_{L^2(\R^2)}\leq 2\mathcal{E}[v_0]$. On the other hand, in the focusing case $\sigma=1$, if we assume that $\norm{v_0}_{L^2(\R^2)} < \norm{\varphi}_{L^2(\R^2)}$, then the sharp Gagliardo-Nirenberg inequality \eqref{GNcrit} yields
$$
\| {\nabla v(t)}\|^2_{L^2(\R^2)}\leq \frac{2}{\left(1-\frac{\norm{v_0}_{L^2(\R^2)} ^2}{ \norm{\varphi}_{L^2(\R^2)}^2}\right)}\mathcal{E}[v_0].
$$
Therefore, the solution of \eqref{agzk} is global in $H^1(\R^2)$ for all $v_0 \in H^1(\R^2)$ when $\sigma = -1$, and for all $v_0 \in H^1(\R^2)$ such that $\norm{v_0}_{L^2(\R^2)} < \norm{\varphi}_{L^2(\R^2)}$ when $\sigma = 1$.



In this paper we are interested in global well-posedness question in $H^s(\R^2)$ with $s< 1$, and in particular, in the application of the almost conservation method to the setting of Zakharov-Kuznetsov model. 
We use the $I$-method introduced by Colliander, Keel, Staffilani,
Takaoka and Tao \cite{CKSTT} (see also \cite{CKSTTKdV}). Briefly recalling the approach, we note that since the energy is not well-defined for initial data in $H^s(\R^2)$ with $s< 1$, a smoothing operator $I$ is introduced, so that the energy of the smoothened solution (or, the modified energy) is finite. Even though the modified energy of the solution is not conserved in time, it can be shown to be slowly growing (or, almost conserved). By controlling the growth of this modified energy, we can iterate a local existence result finitely many times to obtain the existence of the solution for any time $T>0$.

One of the main ingredients of proving the type of the local
existence theorem we need is to establish a trilinear estimate 
in the suitable Bourgain spaces $X_{s,b}$ 
associated to the linear part of the equation. To prove
the trilinear estimate, we use a smoothing effect proved by Faminskii
\cite{Faminskii}, which is the 2d upgrade of the smoothing attained
by Kenig, Ponce and Vega \cite{KPVKdV} for the Airy equation. To complete the proof, we need an
$L^4$-based maximal in time estimate, which can be easily obtained by
following Kenig, Ponce and Vega's approach \cite{KPVKdV} for the KdV equation, if the dispersion relation associated to the
underlying Bourgain spaces is symmetric in both spatial variables
(see, for example, Gr\"unrock \cite{Grun2d3dmZK}). Therefore,
following Gr\"unrock and Herr \cite{GH}, we symmetrize the mZK equation and work with the symmetrized version instead. We prove the trilinear estimate in the Bourgain spaces associated with the linear part of the symmetrized mZK equation. Since the symmetrization changes the conserved
quantities associated with the IVP, we also prove that the $H^1(\R^2)$ norm of the smoothened solution is bounded by its energy.
Finally, we obtain a polynomial growth of the $H^s(\R^2)$ norm of the solution.

In this paper, we consider both focusing and defocusing cases of the mZK equation. In the focusing case, $\sigma = 1$, we prove global well-posedness in $H^s(\R^2)$, $s >\frac{3}{4}$, under the assumption that $L^2$ norm of the initial data is less than the $L^2$ norm of the ground state solution. In the defocusing case, $\sigma = -1$, we prove the same result without any restriction on the size of initial data. As discussed before, this is exactly the situation in $H^1(\R^2)$.
Our first result is the following 
\begin{thm}\label{thm:mainOriginal}
The initial value problem (\ref{agzk}) is globally well-posed in $H^s(\R^2)$, $s >\frac{3}{4}$, for any $v_0 \in H^s(\R^2)$ when $\sigma = -1$, and for $v_0 \in H^s(\R^2)$ such that  
\begin{equation}\label{eq:smallnessOriginal}
\norm{v_0}_{L^2(\R^2)} < \norm{\varphi}_{L^2(\R^2)}
\end{equation}
when $\sigma = 1$; here, $\varphi$ is the ground state solution of \eqref{GSeq}.

Furthermore, for any time $T>0$, the solution $v(t)$ satisfies the following polynomial bound
\begin{equation}\label{apolybound}
\sup_{t \in [0,T]} \norm{v(t)}_{H^s (\R^2)} \lesssim \left( 1 + T \right)^{\frac{1-s}{4s -3}+}.
\end{equation}
\end{thm}

To prove the above result we deal with the symmetrized version of the mZK equation. For that we make a linear change of variables $x \mapsto a x +b y$ and 
$y \mapsto a x -b y$ with $a = 2^{-\frac{2}{3}}$ and $b = 3^{\half}2^{-\frac{2}{3}}$, following Gr\"unrock and Herr \cite{GH}, to obtain
\begin{equation}
  \begin{cases}
    \partial_t u + (\partial_{x}^3 + \partial_{y}^3 ) u + \sigma a( \partial_{x} + \partial_{y}) (u^{3}) = 0\\
    u(x,y,0) = u_0,
  \end{cases}
  \label{gzk}
\end{equation}
where $\sigma$ still denotes the sign; see details on this symmetrization and properties of the new equation in Section \ref{S:sym}. Since this change of variables \eqref{cov} is essentially a rotation, we study the equation \eqref{gzk} instead of \eqref{agzk} without changing the well-posedness theory. Our second result is the following 
\begin{thm}
The initial value problem (\ref{gzk}) is globally well-posed in $H^s(\R^2)$, $s >\frac{3}{4}$, for any $u_0 \in H^s(\R^2)$ when $\sigma = -1$, and for  $u_0 \in H^s(\R^2)$ such that  
\begin{equation}
\norm{u_0}_{L^2(\R^2)} < \sqrt{2ab} \, \norm{\varphi}_{L^2(\R^2)},
\label{eq:smallness}
\end{equation}
when $\sigma = 1$. 

Furthermore, for any time $T>0$, the solution $u(t)$ of \eqref{gzk} satisfies the following polynomial bound
\begin{equation}
  \sup_{t \in [0,T]} \norm{u(t)}_{H^s(\R^2)} \lesssim \left( 1 + T
  \right)^{\frac{1-s}{4s -3}+}.
  \label{polybound}
\end{equation}
  \label{thm:main}
\end{thm}

\begin{rem}
Since the Jacobian of the above change of variables 
is $\abs{2ab}$, undoing the change, 
we see that the threshold condition \eqref{eq:smallness} is equivalent to
\eqref{eq:smallnessOriginal}, and \eqref{polybound} implies \eqref{apolybound}. Therefore, Theorem \ref{thm:main} and Theorem \ref{thm:mainOriginal} are equivalent, including the polynomial growth bounds \eqref{polybound} and \eqref{apolybound}.
\end{rem}

This paper is organized as follows. In the next section, we introduce
some notations and preliminaries that will be used throughout the
paper. In Section \ref{trilinear} we prove a trilinear estimate
required for a variant of the usual local existence theorem. Section
\ref{bilinear} contains a refinement of bilinear Strichartz estimate
when frequencies are separated. In Section \ref{conservation}, we
introduce the modified energy and prove the almost conservation law.
Section \ref{variant} is about the local existence theorem of the
modified solution. Finally, in Section \ref{mainproof}, we prove the
global result stated in Theorem \ref{thm:main}, which in turn implies Theorem \ref{thm:mainOriginal}.

\subsection{Acknowledgements}
L.G.F. was partially supported by CNPq, CAPES and FAPEMIG (Brazil). S.R. was partially supported by the NSF CAREER grant DMS-1151618 and NSF grant DMS-1815873.  D.B. had partial graduate research and travel support to work on this project from grant DMS-1151618 (PI: Roudenko).

\section{Notations, symmetrization and linear estimates}
Throughout this paper, we shall denote the two dimensional spatial
variable pair by $(x,y)$ and its dual Fourier variable by $ \zeta = (\xi, \eta)$. For any integer $j$, we denote $\zeta_j = (\xi_j, \eta_j)$.
As usual, we denote the time variable by $t$ and its Fourier dual variable by $\tau$.  

We use $\F_{x,y,t}$ (or $\F$) to denote the Fourier transform both in space and time variables 
\begin{align*}
  \F_{x,y,t}(f) (\xi, \eta ,\tau) = \iiint_{\R^3} f(x,y,t) e^{-i(x \xi + y \eta + t \tau)} dx dy dt.
\end{align*}
When the Fourier transform is computed in one or two variables out of $x$, $y$ and $t$, we write the variable(s) as a subscript of $\F$. 
In this way, for a function $f \in \mathcal{S}$, the Schwartz class, for example, $\F_{x,y}$ is defined by 
\begin{align*}
  \F_{x,y}(f)(\xi, \eta, t) = \iint_{\R^2} f(x,y,t) e^{-i(x \xi + y \eta)} dx dy.
\end{align*}
$\F_{x,y,t}$ and $\F_{x,y}$ are also denoted by $\widehat{(.)}$ and $\widetilde{(.)}$, respectively.

Given $f(\xi,\eta,\tau)$ in the Schwartz class, $\F^{-1}_{\xi, \eta, \tau} (f) (x,y,t)$ denotes the inverse Fourier transform of $f(\xi, \eta, \tau)$ in both space and time variables and the inverse Fourier transform is written in $(x,y,t)$ variables. Similar to $\F$, when the inverse Fourier transform is computed in one or two variables out of $\xi$, $\eta$ and $\tau$, we denote the variables as a subscript of $\F^{-1}$.

By $\norm{.}_{L^p}$ or $\norm{.}_p$ we denote the $L^p(
\R)$ norm. We use a subscript to denote the variable with respect to which norm is computed.  
The mixed Lebesgue norm is defined by
\begin{align*}
  \norm{u}_{L_x^pL_y^qL_t^r} = \left( \int_{\R} \left( \int_{\R} \left(\int_{\R} |u(x,y,t)|^r dt\right)^{\frac{q}{r}} dy \right)^{\frac{p}{q}} dx \right)^{\frac{1}{p}}
\end{align*}
with obvious modifications when $p$, $q$ or $r$ is $\infty$. We abbreviate $L_x^pL_y^pL_t^q$ and $L_x^pL_y^pL_t^p$ by $L_{x,y}^pL_t^q$ and $L_{x,y,t}^p$, respectively.

By $D^\alpha$ and $J^\alpha$, we define the Fourier multipliers with
symbols $|\zeta|^\alpha$ and $\langle \zeta \rangle^{\alpha}$
respectively, where $\langle \zeta \rangle = \sqrt{1 + |\zeta|^2}$.
Also, $D_x^\alpha$ and $D_y^\alpha$ denote the Fourier multipliers
with symbols $|\xi|^\alpha$ and $|\eta|^\alpha$, respectively.
$J_y^\alpha$ and $J^\alpha_y$ denote the Fourier multipliers
with symbols $\jap{\xi}^\alpha$ and $\jap{\eta}^\alpha$, respectively.
In this notation, the norm in the Sobolev spaces $H^s(\R^2)$ and $\dot{H}^s(\R^2)$ are defined by
\begin{align*}
  \norm{u}_{H^s(\R^2)} = \norm{J^s u}_{L^2(\R^2)}
\end{align*}
and 
\begin{align*}
  \norm{u}_{\dot{H}^s(\R^2)} = \norm{D^s u}_{L^2(\R^2)}.
\end{align*}

Let $s,\ b \in \R$. The Bourgain space $X_{s,b}$ is defined as the space of all tempered distributions $u$ on $\R^2\times \R$ such that
\begin{align*}
  \norm{u}_{X_{s,b}} = \norm{ \langle \zeta \rangle^s \langle \tau -
  (\xi^3 + \eta^3) \rangle^b \F(u)(\xi, \eta, \tau)}_{L^2_{\xi,\eta,\tau}} < \infty.
\end{align*}
Also, for $T>0$, we define the localized $X_{s,b}^T$ norm by 
\begin{equation}
  \norm{u}_{X_{s,b}^T} = \inf \left\{ \norm{v}_{X_{s,b}} : v \in X_{s,b} \text{ and }  v(t) =
  u(t) \text{ on } [0,T] \right\}.
  \label{def:rest}
\end{equation}
These spaces were first used by Bourgain \cite{Bourg1, Bourg2} to study the nonlinear Schr\"odinger equation and the KdV equation, respectively.

Given $A,B \ge 0$, we write $A \lesssim B$,  if for some universal constant $K > 2$ we have $ A \le KB$. We write $A \sim B$, if both $ A \lesssim B$ and $ B \lesssim A$ hold. We write $A << B$, if there is a universal constant $K > 2$ such that $ KA < B$.  
For arbitrarily small $\varepsilon >0$, we use $a+$ and $a-$ to
denote $a+\varepsilon$ and $a-\varepsilon$, respectively. By $a++$ and
$a- -$ we denote $a+2\varepsilon$ and $a-2\varepsilon$, respectively.

\subsection{Symmetrization Details} \label{S:sym}

To symmetrize the initial mZK equation (\ref{agzk}), we perform the
linear change of variables with  $a = 2^{-\frac{2}{3}}$ and $b = 3^{\half}2^{-\frac{2}{3}}$
\begin{equation}
\begin{cases}
    x'=a x +b y\\
    y'=a x -b y,
\end{cases}
\label{cov}
\end{equation}
and denote $u(x',y') = v(x,y)$ and $u_0(x',y') = v_0(x,y)$. Then
\begin{equation}
  \begin{split}
  \partial_x v(x,y) = a (\partial_{x'} + \partial_{y'}) u(x',y')\\
  \partial_y v(x,y) = b (\partial_{x'} - \partial_{y'}) u(x',y').
  \end{split}
  \label{differentials}
\end{equation}
Thus, re-defining $x'$ by $x$ and $y'$ by $y$, the IVP \eqref{agzk} becomes
\eqref{gzk}
Note that, from (\ref{differentials}), we can write
\begin{align*}
  \abs{\nabla v}^2 = v_x^2 + v_y^2 & = a^2 \left(u_{x'} + u_{y'} \right)^2 + b^2 \left( u_{x'} -  u_{y'} \right)^2 = (a^2 + b^2) (u_{x'}^2 + u_{y'}^2) + 2(a^2 - b^2) u_{ x'}u_{y'}.
\end{align*}
Thus, using the change of variables \eqref{cov}, we get
\begin{equation}
  \norm{\nabla v}_{L^2(\R^2)}^2 = \frac{a^2 + b^2}{|2ab|}
  \norm{\nabla u}_{L^2(\R^2)}^2 - \frac{2(a^2 - b^2)}{|2ab|} \int_{\R^2}^{} u_{x'}u_{y'} dx'dy'.
  \label{covNabla}
\end{equation}
Writing the energy $\mathcal{E}[v]$ in terms of $u$, we get,
\begin{align*}
  \frac{2}{a^2 + b^2} |2ab| \mathcal{E}[v] & =  \norm{\nabla u
  }_{L^2(\R^2)}^2 +  \frac{2(a^2 - b^2)}{a^2 + b^2}
  \int_{\R^2}^{} u_{x'}u_{y'} dx' dy' -\frac{2\sigma}{4(a^2 +
  b^2)} \norm{u}_{L^4(\R^2)}^4.
\end{align*}
Observe that $(a^2 + b^2) = 2^{\frac{2}{3}}$ and $\frac{2(a^2 - b^2)}{a^2 + b^2 } = -1$. We next define the energy of $u$ by 
\begin{equation}
E[u(t)] = \frac{|2ab|}{a^2 + b^2} \mathcal{E}[v(t)].
\label{energyDef}
\end{equation}
Then
\begin{equation}
E[u(t)] = \half \int_{\R^2} |\nabla u(x,y,t)|^2 dxdy - \half \int_{\R^2} (u_xu_y)(x,y,t) dxdy - \frac{\sigma a}{4} \int_{\R^2} u^{4}(x,y,t) dxdy,
  \label{energy}
\end{equation}
where $a = (a^2 + b^2)^{-1} =  2^{-\frac{2}{3}}$.

\begin{rem}
In view of (\ref{energy}), the symmetrized ZK equation \eqref{gzk} can be written in the Hamiltonian form
\begin{align*}
    u_t =  (\partial_x + \partial_y) E'(u).
\end{align*}
\end{rem}
Since $E[u(t)]$ is conserved in time, so is $\mathcal{E}[v(t)]$. From \eqref{amass} we also have mass conservation for $u(t)$ for all time $t >0$
\begin{equation}
  \norm{u(t)}_{L^2(\R^2)} = \norm{u_0}_{L^2(\R^2)}.
  \label{masscons}
\end{equation}


\subsection{Linear estimates}
We denote the unitary group associated to the linear part of symmetrized equation (\ref{gzk}) by 
\begin{align*}
  U(t) = e^{-t(\partial_x^3 + \partial_y^3)},
\end{align*}
that is, for any $u_0 \in H^s(\R^2)$ the propagator $U(t)u_0$ is the solution to the linear problem
\begin{equation}
\label{gzk-linear}
  \begin{cases}
  \partial_t u + (\partial_{x}^3 + \partial_{y}^3 ) u = 0\\
  u(x,y,0) = u_0(x,y).
  \end{cases}
\end{equation}

Next, we recall some linear estimates in the mixed Lebesgue spaces as
well as in the $X_{s,b}$ spaces. 
Following the proof of Theorem 2.2 of Faminskii \cite{Faminskii} and using the fact that $(\xi^3 + \eta^3)$ is monotonic in both $\xi$ and $\eta$ variables, we get the smoothing estimates
\begin{equation*}
  \norm{D_x U(t)u_0}_{L^\infty_x L^2_{y,t}} \lesssim \norm{u_0}_{L^2(\R^2)}
\end{equation*}
and
\begin{equation}
  \norm{D U(t)u_0}_{L^\infty_x L^2_{y,t}} \lesssim \norm{u_0}_{L^2(\R^2)}.
  \label{smoothing}
\end{equation}

From Theorem 3.1.(ii) of Kenig, Ponce and Vega \cite{KPVOsc}
(see also estimate (7) of Gr\"unrock and Herr \cite{GH}), we have the Strichartz type estimate
\begin{equation}
  \norm{D_x^{\frac{1}{2p}} D_y^{\frac{1}{2p}} U(t)
u_0}_{L^p_tL^q_{x,y}}\le \norm{u_0}_{L^2(\R^2)},
  \label{str1}
\end{equation}
where $\frac{2}{p} + \frac{2}{q} = 1$, $p>2$.
Taking $p=q=4$, we have,
\begin{equation*}
  \norm{D_x^{\frac{1}{8}} D_y^{\frac{1}{8}} U(t)
u_0}_{L^4_{x,y,t}}\le \norm{u_0}_{L^2(\R^2)}.
\end{equation*}

Also, taking $p=5$, $q = 10/3$ in (\ref{str1}), we deduce
\begin{equation}
  \norm{D^{\frac{1}{10}}_x D^{\frac{1}{10}}_y U(t)u_0}_{L_t^5L_{x,y}^{\frac{10}{3}}} \lesssim \norm{u_0}_{L^2(\R^2)}.
  \label{str2}
\end{equation}
From (\ref{str2}), using Sobolev embedding in $x$ and $y$ variables,
and applying Lemma 2.3 of Ginibre, Tsutsumi and Velo \cite{GTV}, we have
\begin{equation}
\norm{u}_{L^5_{t,x,y}} \le \norm{u}_{X_{0, \frac{1}{2}+}}.
\label{strsob}
\end{equation}
Interpolating (\ref{strsob}) with the trivial estimate
$
\norm{u}_{L^2_{x,y,t}} \le \norm{u}_{X_{0,0}}$,
we obtain
\begin{equation}
\norm{u}_{L^4_{x,y,t}} \le \norm{u}_{X_{0, \frac{1}{2}--}}.
\label{interpolated}
\end{equation}

We recall the classical inequality (see Ginibre, Tsutsumi and Velo \cite{GTV} or  Tao's book \cite[Corollary 2.10]{TaoBook})
\begin{equation}
  \norm{u}_{L^{\infty}_tL^{2}_{x,y}} \lesssim \norm{u}_{X_{0, \half+}}.
  \label{classical}
\end{equation}

Applying the Sobolev embedding $H^{1+}(\R^2) \hookrightarrow L^\infty(\R^2)$ to (\ref{classical}), we have
\begin{equation}
  \norm{u}_{L^\infty_{t,x,y}} \lesssim \norm{D^{1+}u}_{X_{0, \half+}}.
  \label{infinf}
\end{equation}

Interpolation between  (\ref{strsob}) and (\ref{infinf}) yields
\begin{equation}
  \norm{u}_{L^p_{t,x,y}} \lesssim \norm{D^{\alpha(p)} u}_{X_{0, \half+}},
  \label{alphap}
\end{equation}
where $\alpha(p) = (1+)\left(\frac{p-5}{p}\right)$ and $p \in (5, \infty)$.

\section{A Trilinear Estimate}
\label{trilinear}
In this section, we prove a trilinear estimate, which will be the
key ingredient in the proof of the local well-posedness theory of the modified solution $I_Nu$ (see Theorem \ref{thm:xsb} below).
First, we state the $L^4$-based maximal in time estimate from 
Gr\"unrock \cite{Grun2d3dmZK}.
\begin{lemma}[Proposition 1 in \cite{Grun2d3dmZK}]
Let $U(t)u_0$ be the free solution to the linear symmetrized ZK equation \eqref{gzk-linear}. Then
\begin{equation}
    \norm{U(t)u_0}_{L_{x,y}^4L_t^\infty} \lesssim
    \norm{D_x^{\frac{1}{4}}D_y^{\frac{1}{4}}u_0}_{L^2(\R^2)}.
    \label{gain14}
\end{equation}
\end{lemma}

An immediate consequence of the previous result is the following lemma.

\begin{lemma}
Let $U(t)u_0$ be the free solution to the linear symmetrized ZK equation \eqref{gzk-linear}. Then
\begin{equation}
\norm{U(t)u_0}_{L_{x}^4L_{y,t}^\infty} \lesssim \norm{D^{\frac{3}{4}+}u_0}_{L^2(\R^2)}.
\label{maximal-4}
\end{equation}
\end{lemma}
\proof
Since $\xi\eta \lesssim |\xi|^2 + |\eta|^2$, we have $|\xi \eta|^{\frac{1}{4}} \lesssim |\zeta|^{\half}$. Hence, \eqref{gain14} implies
\begin{equation}
    \norm{U(t)u_0}_{L_{x,y}^4L_t^\infty} \lesssim \norm{D^{1/2}u_0}_{L^2(\R^2)}. 
    \label{cor1}
\end{equation}
On the other hand, taking $L_x^4L_t^\infty$ in the both sides of
the following Sobolev inequality in $y$ variable (see \cite[Page 336]{TaoBook}) 
\begin{equation*}
  \norm{u}_{L_y^\infty} \lesssim \| {D^{\frac{1}{4}+}_y u} \|_{L_y^4},
\end{equation*}
we get
\begin{equation}
  \norm{U(t)u_0}_{L_x^4 L_t^\infty L_y^\infty} \lesssim
  \| {D_y^{\frac{1}{4}+}U(t)u_0} \|_{L_x^4L_t^\infty L_y^4}.
  \label{sob2}
\end{equation}
Now, using Minkowski's inequality, 
we can interchange the $L_t^\infty$ and $L^4_y$ norms in the right-hand side of (\ref{sob2}) to get
\begin{equation}
  \norm{U(t)u_0}_{L_x^4 L_{y,t}^\infty} \lesssim
  \norm{D_y^{\frac{1}{4}+}U(t)u_0}_{L_{x,y}^4L_t^\infty }.
  \label{sobmin}
\end{equation}
Then,  from (\ref{sobmin}) and  (\ref{cor1}), we have
\begin{equation*}
  \norm{U(t)u_0}_{L_{x}^4 L_y^\infty L_t^\infty} \lesssim
  \norm{D^{\half} D_y^{\frac{1}{4}+}u_0}_{L^2(\R^2)}.
\end{equation*}
Since $\abs{\eta} \le \abs{\zeta}$,  (\ref{maximal-4}) follows.
\qed 


Now, we prove the following trilinear estimate.
\begin{lemma} 
  For any $\frac{3}{4}<s<1$, the following inequality holds
\begin{equation}
  \norm{(\partial_x + \partial_y)(u_1 u_2 u_3) }_{X_{s,-\frac{1}{2}++}} \le \prod\limits_{i=1}^3 \norm{u_i}_{X_{s,\frac{1}{2}+}}.
  \label{basic:trilinear}
  \end{equation}
  \label{lemma:trilinear}
\end{lemma}
\proof
Since $|\xi| + |\eta| \lesssim \sqrt{|\xi|^2 + |\eta|^2}$, it is enough to show that
\begin{align*}
  \norm{D\left( \prod_{j=1}^{3} u_j \right)}_{X_{s,-\frac{1}{2}++}} \lesssim
  \prod_{j=1}^3 \norm{u_j}_{X_{s,\half+}}.
\end{align*}
We proceed as in the work of the second author \cite{FcgKdV}.
First note that it is enough to show
\begin{align*}
  \norm{D\left( \prod_{j=1}^{3} u_j \right)}_{X_{s,0}} \lesssim
  \prod_{j=1}^3 \norm{u_j}_{X_{s,\half+}}.
\end{align*}
Next, from the definition,
\begin{align*}
  \norm{D\left( \prod_{j=1}^{3} u_j \right)}_{X_{s,0}}  =
  \norm{\abs{\zeta} \langle \zeta \rangle^s \int_{\ast}^{}
  \prod_{j=1}^{3}\widetilde{u}_j(\zeta_j, \tau_j) d\zeta_1 d\tau_1
d\zeta_2 d\tau_2}_{L^2_{\xi, \eta, \tau}},
\end{align*}
where $\ast$ denotes the set  $\{\sum_{i=1}^3(\xi_i,\eta_i,\tau_i) = (\xi,\eta,\tau)\}$.
We decompose the domain of integration according to the relative sizes of spacial frequencies. By symmetry, we assume $|\zeta_1| \ge |\zeta_2| \ge |\zeta_3|$.
We consider the following 3 regions
\begin{align*}
   A & = \{ |\zeta_1| \le 1\},\\
   B & = \{ |\zeta_1| \ge 1,\ |\zeta_3| < |\zeta_1|/2\},\\
   C & = \{ |\zeta_1| \ge 1,\ |\zeta_3| \ge |\zeta_1|/2\}.
\end{align*}

In region A, we have 
\begin{align*}
  \abs{\abs{\zeta} \langle \zeta \rangle^s} \le \left(\sum_{i=1}^3
  |\zeta_i| \right) \left(1+\sum_{i=1}^3 |\zeta_i|^2 \right)^\half \lesssim 1.
\end{align*}
Hence, via Plancherel's theorem and Holder's inequality,
\begin{align*} 
\norm{\abs{\zeta} \langle \zeta \rangle^s \int_{*}^{} \prod_{j=1}^{3}   \widetilde{u}_j(\zeta_j, \tau_j)d\zeta_1 d\tau_1 d\zeta_2 d\tau_2}_{L^2_{\xi, \eta, \tau}} 
  & \lesssim \norm{\int_{*}^{} \prod\limits_{j=1}^{3} \widetilde{u}_j(\zeta_j,\tau_j)d\zeta_1 d\tau_1 d\zeta_2 d\tau_2}_{L^2_{\xi, \eta, \tau}} 
   = \norm{\prod\limits_{j=1}^{3} u_j}_{L^2_{x,y,t}} \\
   & \lesssim \prod\limits_{j=1}^{3} \norm{u_j}_{L^6_{x,y,t}}
   \lesssim \prod\limits_{j=1}^{3} \norm{u_j}_{X_{\alpha(6), \half+}} 
   \lesssim \prod\limits_{j=1}^{3} \norm{u_j}_{X_{s,\half+}},
\end{align*}
where we have used \eqref{alphap} with $p = 6$ and $\alpha(6) = \frac{1}{6}+ < 3/4 <s$.

In region B,  it is easy to see that $|\zeta|\langle \zeta \rangle^s \lesssim |\zeta_1| \langle \zeta_1 \rangle^s$. 
Moreover, using Lemma 2.3 in \cite{GTV}, from the smoothing estimate
\eqref{smoothing} and the $L_x^4$-maximal function estimate \eqref{maximal-4}, we get
\begin{equation}
  \norm{D u}_{L^\infty_x L^2_{y,t}} \lesssim \norm{u}_{X_{0, \half+}}
  \label{smoothing-xsb}
\end{equation}
and
\begin{equation}
  \norm{u}_{L_{x}^4L_{y,t}^\infty} \lesssim  \norm{u}_{X_{\frac{3}{4}+, \half+}}.
  \label{maximal-4-xsb}
\end{equation}
Applying H\"older's inequality, \eqref{smoothing-xsb} and \eqref{maximal-4-xsb}, we get
\begin{align*}
  \norm{\abs{\zeta} \langle \zeta \rangle^s \int_{*}^{} \prod\limits_{j=1}^{3} \widetilde{u}_j(\zeta_j, \tau_j)d\zeta_1 d\tau_1 d\zeta_2 d\tau_2}_{L^2_{\xi, \eta, \tau}} & \lesssim \norm{ \int_{*}^{} \abs{\zeta_1} \langle \zeta_1 \rangle^s \prod\limits_{i=1}^{3} \widetilde{u}_j(\zeta_j,\tau_j)d\zeta_1 d\tau_1 d\zeta_2 d\tau_2}_{L^2_{\xi, \eta, \tau}} \\
  & = \norm{ (D J^s u_1) u_2 u_3}_{L^2_{x,y,t}}\\
  & \lesssim \norm{DJ^su_1}_{L_x^\infty L_{y,t}^2} \norm{u_1}_{L_x^4
    L_{y,t}^\infty} \norm{u_2}_{L_x^4 L_{y,t}^\infty}\\
  & \lesssim \norm{u_1}_{X_{s,\half+}} \prod_{i=2}^3
  \norm{u_i}_{X_{\frac{3}{4}+,\half+}}.
\end{align*}

In region C, we have \begin{align*}
  \zeta_1 \sim \zeta_2 \sim \zeta_3.
\end{align*}
Thus,
\begin{align*}
  \abs{\sum \zeta_i} \jap{ \sum \zeta_i }^s  \lesssim  \jap{\zeta_1}^s \langle \zeta_2 \rangle^\half \langle \zeta_3 \rangle^\half.
\end{align*}
Hence, using H\"older's inequality, we have
\begin{align*}
  \norm{\abs{\zeta} \langle \zeta \rangle^s \int_{*}^{} \prod\limits_{j=1}^{3} \widetilde{u}_j(\zeta_j, \tau_j)d\zeta_1 d\tau_1 d\zeta_2 d\tau_2}_{L^2_{\xi, \eta, \tau}} 
  & \lesssim \norm{ \int_{\ast}^{} \langle \zeta_1 \rangle^s \widetilde{u}_1 \langle \zeta_2 \rangle^\half \widetilde{u}_2 \langle \zeta_3 \rangle^\half \widetilde{u}_3 d\zeta_1 d\tau_1 d\zeta_2 d\tau_2}_{L^2_{\xi, \eta, \tau}}\\
  & = \norm{ J^s u_1 J^\half u_2 J^\half u_3 }_{L^2_{x,y,t}}\\
  & \lesssim \norm{J^su_1}_{L^5} \norm{J^{\half}
u_2}_{L^{\frac{20}{3}}} \norm{J^\half u_3}_{L^{\frac{20}{3}}}\\
  & \lesssim \norm{u_1}_{X_{s,\half+}}
  \norm{u_2}_{X_{\frac{3}{4}+,\half+}} \norm{u_3}_{X_{\frac{3}{4}+,\half+}},
\end{align*}
where we have used the $L^5$ Strichartz estimate (\ref{strsob}) and
the inequality \eqref{alphap} with $p= \frac{20}{3}$ (note that $\alpha(\frac{20}{3}) + \frac{1}{2} = (\frac{1}{4}+) + \frac{1}{2} = \frac{3}{4} +$).
\qed

\section{A refinement of bilinear Strichartz estimate }
\label{bilinear}
In this section, we prove a refinement of the bilinear Strichartz estimate when 
supports of frequencies are separated.
Originally, Bourgain \cite{Bourgain} introduced such a refinement 
in the context of the nonlinear Schr\"odinger equation in two space dimensions. 
The second authot \cite{FcgKdV} also derived an improvement of the bilinear Strichartz estimate associated with the KdV equation for frequencies that are separated.
It should be mentioned 
that Molinet and Pilod \cite{MP} proved a similar
refinement for the unitary group associated with the linear part of
the ZK equation (\ref{agzk}) using a dyadic decomposition. Here, we follow a different approach and prove a refinement for the linear part associated with the symmetrized mZK equation (\ref{gzk}).
\begin{lemma}    \label{lemma:bilWoFreq}
Let $u$ and $v$ be supported on frequencies $\zeta_1$ and $\zeta_2$,    respectively, with $\zeta_1 \sim N_1$, $\zeta_2 \sim N_2$ and $N_1 << N_2$.
Then, we have 
\begin{equation}
     \norm{D^{-\half}uDv}_{L^2_{x,y,t}} \lesssim
     \norm{u}_{X_{0,\half+}} \norm{v}_{X_{0,\half+}}.
     \label{bilWoFreq}
\end{equation}
\end{lemma}
\proof Let $u = U(t)u_0$ and $v= U(t)v_0$. It suffices to show that
\begin{align*}
   \norm{D^{-\half}u Dv }_{L^2_{x,y,t}} \lesssim
   \norm{u_0}_{L^2(\R^2)} \norm{v_0}_{L^2(\R^2)}.
\end{align*}
Recall that $\zeta = (\xi, \eta)$ and $\zeta_i = (\xi_i, \eta_i)$ for $i = 1,2$.
We make change of variables $ \zeta_2 := \zeta - \zeta_1, \zeta_1 := \zeta_1$ to write
\begin{align*}
D^{-\half} u Dv (x,y,t)
   & = \F_{\xi,\eta}^{-1} \left( \widehat{(D^{-\frac{1}{2}}u)} * \widehat{(D^{}v)} (\zeta) \right)(x,y) \\
   & = \int_{\R^2}^{} e^{i(x,y).\zeta}\left( \int_{\R^2}^{} |\zeta_1|^{-\frac{1}{2}}\widehat{U(t)u_0}(\zeta_1)|\zeta - \zeta_1| \widehat{U(t) v_0}(\zeta - \zeta_1) d\zeta_1 \right) d\zeta \\
   & = \int_{\R^2}^{} e^{i(x,y).\zeta}\left( \int_{\R^2}^{} |\zeta_1|^{-\frac{1}{2}} \widehat{U(t)u_0}(\zeta_1)|\zeta - \zeta_1| \widehat{U(t) v_0}(\zeta - \zeta_1) d\zeta_1 \right) d\zeta \\
   & = \int_{\R^2}^{} e^{i(x,y).(\zeta_1 + \zeta_2)}\left( \int_{\R^2}^{} |\zeta_1|^{-\frac{1}{2}} \widehat{U(t)u_0}(\zeta_1)|\zeta_2| \widehat{U(t) v_0}(\zeta_2) d\zeta_1 \right) d\zeta_2 \\
   & = \iiiint_{\R^4} |\zeta_1|^{-\half}|\zeta_2| e^{ix(\xi_1+\xi_2)+ iy(\eta_1+\eta_2)+it(\xi_1^3 + \xi_2^3 + \eta_1^3 + \eta_2^3)} \widehat{u_0}(\zeta_1)\widehat{v_0}(\zeta_2) d\xi_1 d\xi_2 d\eta_1 d\eta_2.
\end{align*}
Now, we make the change of variables $\varphi(\xi_1, \eta_1, \xi_2,\eta_2) = (a,b,c,d)$, where
\begin{align*}
   a & = \xi_1 + \xi_2\\
   b & = \eta_1 + \eta_2\\
   \tau & = \xi_1^3 + \xi_2^3 + \eta_1^3 + \eta_2^3 \\
   c & = \xi_1 - \eta_1.
\end{align*}
Then, using  $|\zeta_i|^2 = \xi_i^2 + \eta_i^2$ for $i=1,2$, the Jacobian of the change of variables $\varphi$ is given by
\begin{align*}
J =  \frac{\partial(a,b,\tau,c)}{\partial(\xi_1,\xi_2,\eta_1,\eta_2)} =
   \det 
   \begin{bmatrix}
     1 & 1 & 0 & 0 \\
     0 & 0 & 1 & 1 \\
     3\xi_1^2 & 3\xi_2^2 & 3\eta_1^2 & 3\eta_2^2 \\
     1 & 0 & -1 & 0
   \end{bmatrix}
   = 3 (\xi_1^2 - \xi_2^2 + \eta_1^2 - \eta_2^2) = 3(|\zeta_1|^2 - |\zeta_2|^2).
\end{align*}
Now, since $N_1 <<  N_2$, we have
\begin{equation}\label{star2}
|J| = 3 \left\vert |\zeta_1|^2 - |\zeta_2|^2\right\vert \gtrsim N_2^2 \ge N_1^2.
\end{equation}
Thus, 
\begin{align*}
D^{-\half} u Dv (x,y,t) = \iiiint_{\R^4} |{\zeta_1}|^{-\half} |{\zeta_2}|e^{ixa+ iby + it\tau} \frac{{\widehat{u}_0}{\widehat{v}_0}}{|J|} dc \, da\, db\, d\tau \\
= \F^{-1}_{a,b,\tau} \left(\int_{\R} |{\zeta_1}|^{-\half} |{\zeta_2}| \frac{{\widehat{u}_0}{\widehat{v}_0}}{|J|} dc\right) (x,y,t).
\end{align*}
Using Plancherel's theorem in time and space variables, Cauchy-Schwarz inequality and \eqref{star2}, we have
\begin{align*}
\norm{D^{-\half} uDv }^2_{L^2_{x,y,t}} & =  \iiint_{\R^3}
   \abs{ \int_{\R} |{\zeta_1}|^{-\half} |{\zeta_2}| \frac{{\widehat{u}_0}{\widehat{v}_0}}{|J|} dc }^2 da\, db\, d\tau\\
   & \le \iiint_{\R^3} \left( \int_{\R}
   \frac{\abs{ \zeta_2 }}{\abs{\zeta_1} \abs{J}^{\half}}dc \right) \left(
  \int_{\R}  {\zeta_2} \frac{\abs{
     {\widehat{u}_0}{\widehat{v}_0} }^2}{\abs{J}^{\frac{3}{2}}} dc\right) da \, db\, d\tau.
\end{align*}
Since $ c = \xi_1 - \eta_1$, we have $\abs{c} \le |\xi_1|+ |\eta_1| \lesssim |\zeta_1| \lesssim N_1$. From (\ref{star2}), we obtain, 
\begin{align*}
  \int_{\R} \frac{|\zeta_2|}{|\zeta_1||J|^{\half}}dc  \lesssim \frac{N_2}{N_1 N_2}\int_{|c| \lesssim N_1} dc \lesssim 1.
\end{align*}
Now, changing the variables $a,b,c,d$ back to $\xi_1,\xi_2, \eta_1, \eta_2$ and using (\ref{star2}) again, we deduce,
\begin{align*}
  \norm{D^{-\half} uDv }^2_{L^2_{x,y,t}} & \lesssim \iiiint_{\R^4}  \frac{ |{\zeta_2}| \abs{ {\widehat{u}_0}{\widehat{v}_0} }^2}{|J|^{\frac{3}{2}}}  dc\, da\, db\, d\tau \\
  & = \iiiint_{\R^4}  \frac{|\zeta_2| \abs{ \widehat{u}_0\widehat{v}_0 }^2}{|J|^{\frac{1}{2}}}  d\xi_1 d\xi_2 d\eta_1 d\eta_2 \\
  & = \norm{u_0}_{L^2(\R^2)}^2 \norm{v_0}_{L^2(\R^2)}^2,
\end{align*}
which concludes the proof of the lemma.
\qed

\begin{cor}
Let $u$ and $v$ be as in the hypothesis of Lemma \ref{lemma:bilWoFreq}. Then we have the following bilinear refinement of Strichartz estimate
\begin{equation}
     \norm{uv}_{L_{x,y,t}^2} \lesssim \frac{N_1^\half}{N_2}
     \norm{u}_{X_{0,\half+}} \norm{v}_{X_{0,\half+}}.
     \label{bil}
\end{equation}
\end{cor}
\proof 
Replacing $u$ by $D^{\half}u$ and $v$ by $D^{-1}v$ in (\ref{bilWoFreq}), we have
\begin{align*}
     \norm{uv}_{L_{x,y,t}^2} \lesssim \norm{D^\half u}_{X_{0,\half+}}
     \norm{D^{-1}v}_{X_{0,\half+}} \lesssim \frac{N_1^\half}{N_2}
     \norm{u}_{X_{0,\half+}} \norm{v}_{X_{0,\half+}}.
\end{align*}
\qed

\section{Modified energy and almost conservation law} \label{conservation}
In this section, we define the modified energy $E^1$ for a solution to the IVP (\ref{gzk}) and prove an almost conservation law.  
We denote the hyperplane $\Gamma_n$ of $(\R^2)^n$ by
\begin{align*}
  \Gamma_n = \left\{ \left( \zeta_1, \dots, \zeta_n \right) \in (\R^{2})^n : \zeta_1 + \zeta_2 + \dots + \zeta_n = 0 \right\},
\end{align*}
equipped with the measure 
\begin{align*}
  \int_{\Gamma_n}^{} f := \int_{(\R^{2})^{(n-1)}}^{} f(\zeta_1,  \dots, \zeta_{n-1}, -\zeta_1 -\dots -\zeta_{n-1}) d\zeta_1 \dots d\zeta_{n-1}
\end{align*}
for any measurable function $f:(\R^{2})^n \rightarrow \mathbb{C}$. For a natural number $n$, we define $S_n$ as the group of all permutations on the set $\{1, \dots, n\}$. A measurable function $M_n$ defined on $\Gamma_n$ is called a
symmetric $n$-multiplier if for any permutation $\pi \in S_n$, we have
\begin{align*}
M_n\left( \zeta_1, \dots, \zeta_n \right) = M_n\left( \zeta_{\pi(1)}, \dots, \zeta_{\pi(n)}\right).
\end{align*}

First, we derive a formula for the time-derivative of a general $n$-linear functional.
\begin{propn}\label{propn:1}
Let $u$ solve the IVP (\ref{gzk}) and $M_n$ be a symmetric $n$-multiplier. 

Let $\Lambda_n$ be an $n$-linear functional defined by
\begin{equation}\label{lambdan}
\Lambda_n(M_n) = \int_{\Gamma_n}^{}M_n(\zeta_1,\dots,\zeta_n) \widehat{u}(\zeta_1)\dots\widehat{u}(\zeta_n).
\end{equation}
Then, recalling the notation $\zeta_i = (\xi_i, \eta_i)$, we have
\begin{align}
&\frac{d}{dt}\Lambda_n(M_n)  = \Lambda_n(M_n\alpha_n) \label{diff-formula}\\
 &- \sigma \,n\, a \, i \Lambda_{n+2}( M_n(\zeta_1,\dots, \zeta_{n-1}, \zeta_n+\dots+\zeta_{n+2}) ( \xi_n + \dots + \xi_{n+2}+ \eta_n + \dots+ \eta_{n+2}) ), \nonumber
\end{align}   
where
\begin{equation}
  \alpha_n = i(\xi_1^3 + \eta_1^3 + \dots + \xi_n^3 + \eta_n^3).
  \label{alpha}
\end{equation}
\end{propn}
\proof
Applying Fourier transform to (\ref{gzk})  in spatial variable, we get
\begin{equation}\label{fourier}
\widehat{u}_t(\xi, \eta) = i[\, (\xi^3 + \eta^3) \widehat{u} - \sigma  a (\xi + \eta)\widehat{u^{3}} \,\,].
\end{equation}

Differentiating (\ref{lambdan}) with respect to $t$ and applying (\ref{fourier}), we get
\begin{align*}
  & \frac{d}{dt}\Lambda_n(M_n) \\
  & = \sum_{j=1}^{n}\int_{\Gamma_n}^{}M_n(\zeta_1,\dots,\zeta_n) \uhat_t(\zeta_j) \prod_{k \ne j} \widehat{u}(\zeta_k)\\
  & = \sum_{j=1}^{n} i \int_{\Gamma_n}^{}  M_n(\zeta_1,\dots,\zeta_n)  \left[ (\xi_j^3 + \eta_j^3) \widehat{u}(\zeta_j) - \sigma a (\xi_j+\eta_j) \widehat{u^{3}}(\zeta_j)\right] \prod_{k \ne j } \uhat(\zeta_i) \\
  & = \sum_{j=1}^{n} i \int_{\Gamma_n}^{}  M_n(\zeta_1,\dots,\zeta_n)  (\xi_j^3 + \eta_j^3) \prod_{k=1 }^{n} \uhat(\zeta_k) 
   - \sum_{j=1}^{n} i\sigma a  \int_{\Gamma_n}^{}  M_n(\zeta_1,\dots,\zeta_n)   (\xi_j+\eta_j) \widehat{u^{3}}(\zeta_j) \prod_{k \ne j } \uhat(\zeta_k). 
\end{align*}  
Now applying a permutation of indexes, interchanging $j$ with $n$ (or the variables $\zeta_j$) in the integration  in the last term above, we rewrite the last term and obtain
\begin{align*}
\frac{d}{dt}\Lambda_n(M_n)   = &  \int_{\Gamma_n}^{}  M_n(\zeta_1,\dots,\zeta_n)  \left(i \sum_{j=1}^{n} \xi_j^3 + \eta_j^3\right) \prod_{k=1 }^{n} \uhat(\zeta_k)\\ 
  & - \sum_{j=1}^{n} i\sigma a  \int_{\Gamma_n}^{}  M_n(\zeta_1,\dots,\zeta_n)   (\xi_n+\eta_n) \widehat{u^{3}}(\zeta_n) \prod_{k=1}^{n-1} \uhat(\zeta_k) \\
   = & \Lambda_n(M_n \alpha_n) - \sigma  \,n\, a \, i \int_{\Gamma_n}^{}  M_n(\zeta_1,\dots,\zeta_n)   (\xi_n+\eta_n) \widehat{u^{3}}(\zeta_n) \prod_{k=1}^{n-1} \uhat(\zeta_k) \\
  = & \Lambda_n(M_n \alpha_n) \\
  & - \sigma  \,n\, a \, i \int_{\Gamma_{n+2}}\!\!\!\!\!\!\! M_n(\zeta_1,\dots,\zeta_{n-1}, \zeta_n+ \dots +\zeta_{n+2}) (\sum_{j=n}^{n+2}\xi_j+\eta_j) \prod_{k=1}^{n+2} \uhat(\zeta_k) \\
  = & \Lambda_n(M_n\alpha_n) \\
  & -\sigma  \,n\, a \, i \Lambda_{n+2}\left( M_n\left( \zeta_1,\dots, \zeta_{n-1}, \zeta_n+\dots+\zeta_{n+2}\right) \left( \xi_n + \dots + \xi_{n+2}+ \eta_n + \dots+ \eta_{n+2}\right) \right),
\end{align*}
where $\alpha_n$ is defined in \eqref{alpha}.
\qed

\subsection{Modified energy functional}
For a large positive real number $N$ and $0<s<1$, define a Fourier
multiplier operator $I_N: H^s(\R^2)  \rightarrow H^1(\R^2)$ by
\begin{equation*}
  \widehat{I_Nf}(\zeta) = m_N(\zeta) \widehat{f}(\zeta),
\end{equation*}
where $m_N$ is smooth, radially symmetric, non-increasing function of $|\zeta|$ such that
\begin{equation*}
  m_N(\zeta) = 
  \begin{cases}
    1 & \text{ if } |\zeta| \le N\\
    \left( \frac{N}{|\zeta|} \right)^{1-s} & \text{ if } |\zeta| \ge 2N .
  \end{cases}
\end{equation*}
Note that for any $f \in X_{s_0, b_0}$, we have
\begin{align}\label{Ismoothing}
  \norm{f}_{X_{s_0, b_0}} \lesssim \norm{I_Nf}_{X_{s_0 +1 -s, b_0}}
  \lesssim N^{1-s} \norm{f}_{X_{s_0, b_0}}.
\end{align}
For simplicity, we drop the subscript $N$ of $m_N$ and only write $m$. Further, for any $j \in \N$, we define
\begin{equation*}
  m_j = m(\zeta_j).
\end{equation*}
Moreover, we sometimes abuse the notation to define $m$ as a function
of a scalar variable $r \in[0,\infty)$ by $m(r) = m(\zeta)$, where $r =|\zeta|$.

We define, for all time $t >0$, the \textit{modified energy} $E^1[u]$ of $u$ by
\begin{equation}\label{star3}
E^1[u](t)= E[I_Nu](t).
\end{equation}

\begin{propn}\label{propn:I_Nugrowth} 
Let $s > \frac{3}{4}$, $N >> 1$ and $u \in H^s(\R^2)$ be a solution to \eqref{gzk} on $[0, \delta]$. Then we have the following growth of the modified energy
\begin{equation*}
\abs{E^1[u](\delta) - E^1[u](0) } \lesssim N^{-1+}  \left( \norm{I_Nu}_{X^\delta_{1,\half+}}^4 +  \norm{I_Nu}_{X^\delta_{1,\half+}}^6 \right).
\end{equation*}
\end{propn}
\proof
First, using Proposition \ref{propn:1}, we write
$\frac{d}{dt}E^1[u](t)$ as a sum of $n$-linear functionals. 
By definition, on $\Gamma_2$, we have, $\zeta_1 = -\zeta_2$ and $\eta_2 = -\eta_1$. Hence, using Plancherel Theorem, we can write
\begin{align*}
  \int_{\R^2}^{}|\nabla u|^2 dxdy 
  &= -\int_{\Gamma_2}^{} \zeta_1 . \zeta_2\  \uhat(\zeta_1) \uhat(\zeta_2)
  = \int_{\Gamma_2}^{}|\zeta_1|^2 \uhat(\zeta_1) \uhat (\zeta_2) \\
  &= \frac{1}{2} \int_{\Gamma_2}^{} (|\zeta_1|^2 + |\zeta_2|^2) \uhat(\zeta_1) \uhat(\zeta_2)
  = \frac{1}{2} \Lambda_2( |\zeta_1|^2 + |\zeta_2|^2),
\end{align*}
and
\begin{align*}
  \int_{\R^2}^{} u_x u_y dxdy 
   &= -\int_{\Gamma_2}^{} \xi_1\eta_2 \uhat(\zeta_1) \uhat(\zeta_2)
   = \int_{\Gamma_2}^{}\xi_1 \eta_1 \uhat(\zeta_1) \uhat(\zeta_2)\\
   &= \frac{1}{2} \int_{\Gamma_2}^{} (\xi_1\eta_1 + \xi_2\eta_2) \uhat(\zeta_1) \uhat(\zeta_2)
   = \frac{1}{2} \Lambda_2(\xi_1\eta_1 + \xi_2\eta_2).
\end{align*}
Finally, one can see that
\begin{align*}
  \int_{\R^2}^{} u^{4} dxdy= \Lambda_{4}(1).
\end{align*}
Therefore, the energy (\ref{energy}) can be written as
\begin{equation}
  E[u] = \frac{1}{4}\Lambda_2(|\zeta_1|^2 + |\zeta_2|^2 - \xi_1\eta_1 - \xi_2\eta_2) - \frac{\sigma  a}{4}\Lambda_{4}(1).
  \label{energyalt}
\end{equation}
Also, using (\ref{energyalt}) and the definition (\ref{star3}), we obtain
\begin{equation} 
  E^1[u]= E[I_Nu] =  \frac{1}{4}\Lambda_2 \left(m_1m_2(|\zeta_1|^2 + |\zeta_2|^2 - \xi_1\eta_1 - \xi_2\eta_2)\right) - \frac{\sigma  a}{4}\Lambda_{4}\left(m_1\dots m_{4}\right).
  \label{menergy}
\end{equation}
Moreover, since $\alpha_2 = 0$ on $\Gamma_2$ and
$\zeta_2+\zeta_3+\zeta_{4} = -\zeta_1$ on $\Gamma_{4}$,
differentiating the first term of (\ref{menergy}) in time and using
\eqref{diff-formula}, we have
\begin{align*}
&  \frac{d}{dt}\Lambda_2\left(m_1\,m_2\left(|\zeta_1|^2 + |\zeta_2|^2 -
\xi_1\eta_1 - \xi_2\eta_2\right)\right) 
= \Lambda_2\left(m_1\,m_2(|\zeta_1|^2 + |\zeta_2|^2 - \xi_1\eta_1 - \xi_2\eta_2) \alpha_2\right) \\
&- 2 i \sigma  a \,\Lambda_{4} \bigg( m_1 \,m( \zeta_2 +\dots+\zeta_{4})
\left( |\zeta_1|^2 + |\zeta_2+\dots+\zeta_{4}|^2 - \xi_1 \eta_1 - (\xi_2 + \dots + \xi_{4})(\eta_2 + \dots + \eta_{4}) \right) \\
&\quad  \times (\xi_2 +\dots+\xi_{4}+\eta_2+\dots+\eta_{4}) \bigg)\\
&  = - 2i\sigma a\Lambda_{4}\left(m_1m(-\zeta_1)\left( |\zeta_1|^2+|-\zeta_1|^2 - \xi_1\eta_1 - (-\xi_1)(-\eta_1) \right)(-\xi_1-\eta_1)\right)\\
&  = 2i\sigma a\Lambda_{4}\left(m_1^2\left( 2|\zeta_1|^2 - 2\xi_1\eta_1 \right)(\xi_1+\eta_1)\right)\\
&  = 4i\sigma a\Lambda_{4}\left(m_1^2(|\zeta_1|^2 - \xi_1\eta_1)(\xi_1+\eta_1)\right)\\
&  = 4i\sigma a\Lambda_{4}\left(m_1^2(\xi_1^3 + \eta_1^3)\right).
\end{align*}
Similarly, differentiating the second term of (\ref{menergy}), we get
\begin{align*}
\frac{d}{dt} &\Lambda_{4}(m_1\dots m_{4})
 = \Lambda_{4}(m_1\dots m_{4} \alpha_{4}) 
   - 4\sigma a i \Lambda_{6}\left(m_1m_2 m_{3} m(\zeta_{4}+\zeta_5+ \zeta_{6})[\xi_{4}+\xi_5 + \xi_{6} + \eta_{4}+ \eta_5 + \eta_{6}]\right) \\
  & = 4i \Lambda_{4}\left(m_1\dots m_{4} (\xi_1^3 + \eta_1^3)\right)
   - 4\sigma a i \Lambda_{6} \left( m_1m_2m_{3} m(\zeta_{4}+\zeta_5 + \zeta_{6})[\xi_{4}+\xi_5 + \xi_{6} + \eta_{4}+ \eta_5 + \eta_{6}] \right),
\end{align*}
where we have used $\Lambda_4(m_1\dots m_4 (\xi_j^3 + \eta_j^3)) =   \Lambda_4(m_1\dots m_4 (\xi_1^3 + \eta_1^3))$ for all $j=1,\dots,4$, by permuting the variables.

Collecting all the  $\Lambda_{4}$ terms appearing in $\frac{d}{dt}E^1[u](t)$, we get
\begin{align*}
      \sigma a i \Lambda_{4}((m_1^2  - m_1\dots m_{4}) (\xi_1^3 +
      \eta_1^3) ) 
      & = \sigma a i \int_{\Gamma_{4}}^{} \frac{m_1^2  - m_1\dots m_{4}}{m_1 \dots m_{4}} (\xi_1^3 + \eta_1^3) \widehat{I_Nu}(\zeta_1) \dots \widehat{I_Nu}(\zeta_{4}) \\
  & = \sigma a i \int_{\Gamma_{4}}^{} \left[  \frac{m_1}{m_2 m_3
  m_{4}} -1 \right] (\xi_1^3 + \eta_1^3) \widehat{I_Nu}(\zeta_1)
  \dots \widehat{I_Nu}(\zeta_{4}).
\end{align*}
Since $\sigma = \pm1$, collecting all the $\Lambda_{6}$ terms, we get
\begin{equation*}
  \begin{split}
  &  \sigma^2  a^2 i \Lambda_{6}(m_1m_2 m_{3} m(\zeta_{4}+ \xi_5 +
  \zeta_{6})[\xi_{4}+ \xi_5 + \xi_{6} + \eta_{4}+ \xi_5 +
  \eta_{6}])\\
  &=    a^2 i \Lambda_{6}(m_1m_2 m_{3} m(\zeta_{4}+\zeta_5+\zeta_{6})[(\xi_{4}+\eta_{4}) + (\xi_5+\eta_5) +  (\xi_{6} + \eta_{6})]).
       \end{split}
\end{equation*}
Hence, the derivative of the modified energy $E^1[u](t)$ is
\begin{align}\label{derivativeEnergy}
       \begin{split}
       \frac{d}{dt} E^1[u](t) & = \sigma a i \int_{\Gamma_{4}}^{} \left[  \frac{m_1}{m_2 m_3 m_{4}} -1 \right] (\xi_1^3 + \eta_1^3) \widehat{I_Nu}(\zeta_1) \dots \widehat{I_Nu}(\zeta_{4}) \\
       & \quad +    a^2 i \Lambda_{6}(m_1m_2 m_{3} m(\zeta_{4}+ \zeta_5 + \zeta_{6}) [(\xi_{4}+\eta_{4}) + (\xi_5+\eta_5) +  (\xi_{6} + \eta_{6})]).
       \end{split}
\end{align}
The Fundamental Theorem of Calculus yields
\begin{align*}
E^1[u](\delta)-E^1[u](0) = \int_{0}^{\delta} \frac{d}{dt} E^1[u](s) \,ds.
\end{align*}
Now, we integrate \eqref{derivativeEnergy} in time variable from $0$ to $\delta$ and take the absolute value to get
\begin{align}\label{terms}
&\abs{E^1[u](\delta) - E^1[u](0)}  \lesssim
\abs{\int_{0}^{\delta} \int_{\Gamma_4}  
\left[ \frac{m_1}{m_2 \,m_3\, m_{4} } -1 \right] (\xi_1^3 + \eta_1^3)
  \prod\limits_{i=1}^{4} \widehat{I_Nu}(\zeta_i,s) ds } \\
& \qquad \qquad + \abs{\int_{0}^{\delta} 
\Lambda_{6} \bigg( m_1\, m_2\dots m_{3} \, m(\zeta_{4}+ \dots + \zeta_{6}) \left[\xi_{4}+\eta_{4}+  \xi_5+\eta_5 +  \xi_{6} + \eta_{6}\right] \bigg) ds}. \nonumber
\end{align}
To estimate the first term on the right-hand side, we decompose the function $u$ into dyadic constituents 
and work with a typical term in the infinite sum. 

For example, we decompose $\widehat{I_Nu}(\zeta_1,t)$ as 
\begin{align*}
       \widehat{I_Nu}(\zeta_1, t) = m(\zeta_1) \widehat{u} (\zeta_1,
       t) =  m(\zeta_1) \sum\limits_{l_1 \in \mathbb{N}}
       \widehat{u_{l_1}}(\zeta_1,t)  = \sum\limits_{l_1 \in \mathbb{N}}
       \widehat{I_Nu_{l_1}}(\zeta_1, t),
\end{align*}
where $\mbox{supp}\,(\widehat{u}_{l_1}) \subset [2^{l_1}, 2^{l_1+1}]$ for each
$l_1 \in \N$. We do a similar decomposition for other functions as well and index the projections by $l_2, l_3, l_4$. 
Then, the first term on the right-hand side of (\ref{terms}) can be written as
\begin{equation}\label{star4}
\sum\limits_{l_1,l_2,l_3,l_4 \in \mathbb{N}} T_{l_1,l_2,l_3,l_4},
\end{equation}
where 
\begin{align*}
T_{l_1, l_2, l_3, l_4} := \int_{0}^{\delta} \int_{\Gamma_{4}}^{} \left[  \frac{m_1}{m_2 m_3 m_{4}} -1 \right] (\xi_1^3 + \eta_1^3) \widehat{I_Nu_{l_1}}(\zeta_1,s )  \widehat{I_Nu_{l_2}}(\zeta_{2},s) \widehat{I_Nu_{l_3}}(\zeta_{3},s) \widehat{I_Nu_{l_4}}(\zeta_{3},s)\,ds
\end{align*}
and $\mbox{supp}\,(\widehat{u}_{l_i}) \subset [N_{l_i}, 2N_{l_i}]$,
$N_{l_i} = 2^{l_i}$ for each $l_i \in \N,\ i = 1,2,3,4$.

Our aim is to show that there exists $\varepsilon >0$ such that
     \begin{align*}
       \abs{T_{l_1,l_2,l_3,l_4}}
       \lesssim N^{-1+} (N_{l_1} N_{l_2} N_{l_3}
       N_{l_4})^{-\varepsilon} \prod\limits_{i=1}^{4}
       \norm{I_Nu_{l_i}}_{X^\delta_{1,\half+}},
     \end{align*}
     so that, after applying the infinite sum over $l_1, \dots, l_4$, we get $N^{-1+}$ on the right-hand side, since $\sum\limits_{l_i \in \N}^{} 2^{-\varepsilon l_i} < \infty$, for $i=1,\dots, 4$.

     Define $N_{max} = \max\left\{N_{l_i}: i =1,2,3,4\right\}$. 
     Since $N_{max}^\varepsilon \ge N_{l_1}^{\frac{\varepsilon}{4}}N_{l_2}^{\frac{\varepsilon}{4}}N_{l_3}^{\frac{\varepsilon}{4}} N_{l_4}^{\frac{\varepsilon}{4}} $, it suffices to show that
     \begin{align*}
       \abs{T_{l_1, l_2,l_3,l_4}} \lesssim N^{-1+} N_{max}^{-\varepsilon} \prod\limits_{i=1}^{4} \norm{I_Nu_{l_i}}_{X^\delta_{1,\half+}}.
     \end{align*}

Here, for brevity we write $u_i$ instead of $u_{l_i}$ for $i=1,2,3,4$. Thus, a typical term in the sum (\ref{star4}) is given by
\begin{equation} \label{Term1}
       Term1 = \abs{\int_{0}^{\delta}
       \int_{\Gamma_{4}}^{} \left[ \frac{m_1}{m_2 m_3 m_{4}}-1
       \right] (\xi_1^3 + \eta_1^3) \widehat{I_Nu_1}(\zeta_1, s) \dots
       \widehat{I_Nu_{4}}(\zeta_{4},s) ds }.
\end{equation}

Without loss of generality, we assume that the Fourier transform of all the functions is non-negative and that $N_2 \ge N_3 \ge N_4$ by the symmetry of the term 
$$ 
\left( \frac{m(\zeta_2+ \zeta_3 + \zeta_4)}{m_2 m_3 m_4} -1 \right) \bigg(-(\xi_2 + \xi_3 + \xi_4)^3 - (\eta_2 + \eta_3 + \eta_4)^3\bigg) 
$$
in $\zeta_2,\zeta_3, \zeta_4$ variables. Moreover, we can assume that $N_2 \gtrsim N$, otherwise, the multiplier is zero, since $m_1 = \dots = m_4 = 1$. Also, the condition $\sum\limits_{i=1}^{4} \zeta_i = 0$ implies $N_1 \lesssim N_2$. \\

       Now, we consider the following nested subcases. Recall that
       our largest frequency is $N_2$. Based on where the second largest frequency $N_3$ is located compared to $N$, we get the following cases
       \begin{enumerate}
	 \item $N_3 << N$. In this subcase, we have $N_2 \gtrsim N >>N_3 \ge N_4$. We also have that $N_1 \sim N_2$. 
	 \item $N \lesssim N_3$. 
	   In this case, based on whether $N_3$ is comparable to $N_2$ or not, we consider the subcases
	   \begin{enumerate}
	     \item $N_3 << N_2 $. This is the case where $ N_{4} \le
	       N_3,\ N \lesssim N_3 << N_2$. Here, we also have $N_1 \sim N_2$.
	     \item $N_3 \sim N_2 $. 
	       This subcase can be broken into two further subcases based on the comparison between  the lowest frequency $N_4$ and $N$
	       \begin{enumerate}
		 \item $ N_4 << N$. Based on the 
		   comparison between $N_1$ and $N_2$, we have two further subcases
		   \begin{enumerate}
		     \item $N_1 << N_2$
		     \item $N_1 \sim N_2$
		   \end{enumerate}
		 \item $N \lesssim N_4$. Again, comparing  $N_4$ and $N_3$, we have the further subcases
		   \begin{enumerate}
		     \item $N_4 << N_3$
		     \item $N_4 \sim N_3$
		   \end{enumerate}
	       \end{enumerate}
	   \end{enumerate}
       \end{enumerate}


Now, we bound the multiplier in every terminal subcase in a pointwise manner.\\

     \textbf{Case 1.} $N_2 \gtrsim N >>N_3 \ge N_{4}$. This implies $N_1 \sim N_2$. By the Mean Value Theorem,  
     \begin{equation*}
       \abs{\frac{m_1}{m_2 m_3 m_4} -1} =  \abs{\frac{m(\zeta_2) -
	 m(\zeta_2+\zeta_3 + \zeta_{4})}{m(\zeta_2)}} \lesssim
	 \frac{N_3}{N_2},
     \end{equation*}
     which we substitute in \eqref{Term1}.
     Using the bilinear Strichartz estimate (\ref{bil}),  we get
\begin{align*}
  Term1 & \lesssim N_1^3 \frac{N_3}{N_2}\norm{I_Nu_1 I_Nu_3}_{L^2(\R^2
  \times [0,\delta])} \norm{I_Nu_2 I_Nu_4}_{L^2(\R^2 \times [0,\delta])} \\
  & \lesssim \frac{N_1^3 N_3}{N_2}
  \frac{1}{\jap{N_3}^{\half}N_1^2}\frac{1}{\langle N_4
    \rangle^{\half}N_2^2} \prod\limits_{i=1}^4 \norm{I_Nu_i}_{X^\delta_{1,\half+}}  \\
  & \lesssim \frac{\jap{N_3}^{\half}}{N_2^2\langle N_4
    \rangle^{\half}} \prod\limits_{i=1}^4 \norm{I_Nu_i}_{X^\delta_{1,\half+}} \text{ (since } N_1 \sim N_2 )\\
  & \lesssim  \frac{N^{\half}}{N_2^{2-} N_2^{0+}}
  \prod\limits_{i=1}^4 \norm{I_Nu_i}_{X^\delta_{1,\half+}} \text{ (since } N_3 \le N \text{ and }  \jap{N_4}^{\half} \ge 1)\\
  & \lesssim  \frac{N^{\half}}{N^{2-} N_{max}^{0+}}
  \prod\limits_{i=1}^4 \norm{I_Nu_i}_{X^\delta_{1,\half+}} \text{ (since } N_2 \gtrsim N)\\
  & \lesssim  N^{-\frac{3}{2}+}N_{max}^{0-} \prod\limits_{i=1}^4
  \norm{I_Nu_i}_{X^\delta_{1,\half+}}.
\end{align*} 

     \textbf{Case 2.(a).} $N_2 >> N_3  \gtrsim N$ and $N_3 \ge N_4$. In this case, we still have $N_1 \sim N_2$.
     
     Here, we use the pointwise bound
     \begin{align*}
\abs{\frac{m_1}{m_2 m_3 m_4} -1} = \abs{\frac{m(\zeta_2+ \zeta_3 + \zeta_{4})}{m(\zeta_2) m(\zeta_3) m(\zeta_{4})} -1} \lesssim \frac{m(\zeta_1)}{m(\zeta_2) m(\zeta_3) m(\zeta_{4})}.
     \end{align*}
     Since $m$ is a non-decreasing function and $N_1 \sim N_2$, we have $m(\zeta_1) \lesssim m(\zeta_2)$. 
     Then, the bound on the multiplier becomes 
     \begin{align*}
       \abs{\frac{m(\zeta_2+ \zeta_3 + \zeta_{4})}{m(\zeta_2) m(\zeta_3) m(\zeta_{4})} -1 }\lesssim \frac{1}{m(N_3) m(N_{4})}.
    \end{align*}

     Hence, using the bilinear Strichartz estimate (\ref{bil}), we have,

     \begin{align*}
       Term1 & \lesssim \frac{N_1^3}{m(N_3)m(N_4)} \norm{I_Nu_1
       I_Nu_3}_{L^2(\R^2 \times [0,\delta])} \norm{I_Nu_2 I_Nu_4}_{L^2(\R^2 \times [0,\delta])} \\
       & \lesssim \frac{N_1^3}{m(N_3)m(N_4)}
       \frac{1}{N_3^{\half}N_1^2} \frac{1}{\jap{N_4}^{\half}
     N_2^2}\prod\limits_{i=1}^4 \norm{I_Nu_i}_{X^\delta_{1,\half+}}\\
       & \lesssim
       \frac{1}{m(N_3)N_3^{\half}m(N_4)\jap{N_4}^{\half}N_2}
       \prod\limits_{i=1}^4 \norm{I_Nu_i}_{X^\delta_{1,\half+}} \text{ (since } N_1 \sim N_2) \\
       & \lesssim \frac{1}{N^{\half}N_2^{1-}N_2^{0+}}
       \prod\limits_{i=1}^4 \norm{I_Nu_i}_{X^\delta_{1,\half+}} \text{ (since } m(N_3)N_3^{\half} \gtrsim N^{\half} \text{ and } m(N_4)\jap{N_4}^{\half} \gtrsim 1) \\
       & \lesssim \frac{1}{N^{\frac{3}{2}-}N_{max}^{0+}}
       \prod\limits_{i=1}^4 \norm{I_Nu_i}_{X^\delta_{1,\half+}} \text{ (since }  N_2 >> N )\\
       & \lesssim  N^{-\frac{3}{2}+}N_{max}^{0-} \prod\limits_{i=1}^4
       \norm{I_Nu_i}_{X^\delta_{1,\half+}} ,
     \end{align*}
     where we have used that for any $p>0$ with $p+s > 1$, the function $m(\zeta)\abs{\zeta}^p$ is increasing and $m(\zeta)\jap{\zeta}^p$ is  bounded below.\\

     \textbf{Case 2.(b).i.A.} In this case, we have $N_4 << N \lesssim N_3 \sim N_2$  and $N_1 << N_2$.\\
     Here, we pair $I_Nu_1$ with $I_Nu_2$ and $I_Nu_3$ with $I_Nu_4$. In the following subcases, we use the following crude bound for the multiplier
     \begin{align*}
       \abs{\frac{m_1}{m_2 m_3 m_4} - 1 }= \abs{\frac{m_1 - m_2 m_3 m_4}{m_2 m_3 m_4}}  \lesssim \frac{1}{m(N_2)m(N_3) m(N_4)}.
     \end{align*}
Thus,
     \begin{align*}
       Term1 & \lesssim \frac{N_1^3}{m(N_2)m(N_3)m(N_4)} \norm{I_Nu_1
       I_Nu_2}_{L^2(\R^2 \times [0,\delta])} \norm{I_Nu_4 I_Nu_3}_{L^2(\R^2 \times [0,\delta])}\\
       & \lesssim \frac{N_1^3}{m(N_2)m(N_3)m(N_4)}
       \frac{N_1^\half}{N_2 \langle N_1 \rangle \langle N_2 \rangle}
       \frac{N_4^\half}{N_3 \langle N_3 \rangle \langle N_4 \rangle}
       \prod\limits_{i=1}^{4} \norm{I_Nu_i}_{X^\delta_{1,\half+}}\\
       & \lesssim \frac{N_1^{\frac{5}{2}}}{m(N_2)m(N_3)m(N_4) N_2^2
	 N_3^2 \jap{N_4}^{\frac{1}{2}}} \prod\limits_{i=1}^{4} \norm{I_Nu_i}_{X^\delta_{1,\half+}}\\
       & \lesssim \frac{1}{m(N_2) N_2 m(N_3) N_3^{\half} m(N_4)
       \jap{N_4}^{\frac{1}{2}}} \prod\limits_{i=1}^{4} \norm{I_Nu_i}_{X^\delta_{1,\half+}}
       \text{ (since } N_1 << N_2 \sim N_3 )\\
       & \lesssim N^{-\frac{3}{2}+} N_{max}^{0-}
       \prod\limits_{i=1}^{4} \norm{I_Nu_i}_{X^\delta_{1,\half+}},
     \end{align*}
     where in the last step we have used $N_{max} \sim N_2$,
     $m(\zeta)\abs{\zeta}^p$ is non-decreasing and
     $m(\zeta)\jap{\zeta}^p > 1$ for $s+p \ge 1$, hence,   $ m(N_3)N_3^\half \gtrsim N^\half,\ m(N_2)N_2 \gtrsim N$ and $m(N_4)\jap{N_4} \gtrsim 1 $.\\

     \textbf{Case 2.(b).i.B.} In this case, we have $N_4 << N \lesssim N_3 \sim N_2 \sim N_1$.\\
     Here, we use the same bound for the multiplier as the previous case.  The only frequencies that are separated are $N_4$ and $N_2$ (or, $N_3$ since $N_2 \sim N_3$). So, we use inequality (\ref{bil})  on the pair $I_Nu_2 I_Nu_4$ and inequality  (\ref{interpolated}) on the other pair to get
     \begin{align*}
       Term1 & \lesssim \frac{N_1^3}{m(N_2) m(N_3) m(N_4)} \norm{I_Nu_2
       I_Nu_4}_{L^2(\R^2 \times [0,\delta])} \norm{I_Nu_3}_{L^4(\R^2
       \times [0,\delta])} \norm{I_Nu_1}_{L^4(\R^2 \times [0,\delta])}\\
       & \lesssim \frac{N_1^3}{m(N_2) m(N_3) m(N_4)}
       \frac{N_4^\half}{N_2 \jap{N_2}\jap{N_4}} \frac{1}{\jap{N_3}
       \jap{N_1}} \prod\limits_{i=1}^{4}\norm{I_Nu_i}_{X^\delta_{1,\half+}}\\
       & \lesssim \frac{N_1^2}{m(N_2) m(N_3) m(N_4) N_2^2
	 \jap{N_4}^\half N_3} \prod\limits_{i=1}^{4}\norm{I_Nu_i}_{X^\delta_{1,\half+}}\\
       & \lesssim \frac{1}{m(N_2) m(N_3) m(N_4) N_2^\half N_3^\half
	 \jap{N_4}^\half} \prod\limits_{i=1}^{4}\norm{I_Nu_i}_{X^\delta_{1,\half+}} 
       \text{ (since } N_1 \lesssim N_2 \sim N_3 \implies N_1^2  \lesssim N_2^\half N_3^\half N_2)\\
       & \lesssim N^{-1+}N_{max}^{0-}\prod\limits_{i=1}^{4}\norm{I_Nu_i}_{X^\delta_{1,\half+}}.
     \end{align*}
     In the last step, we have used $m(N_i)N_i^{\half-} \gtrsim N^{\half-}$,  since $N_i \gtrsim N$ for $i=2,3$, $m(N_4) \jap{N_4}^\half \gtrsim 1$ and $N_{max} \sim N_2$.\\

     \textbf{Case 2.(b).ii.A.} In this case, we have $N \lesssim N_4 << N_3 \sim N_2$.\\
     Here, we again use the same bound for the multiplier as in Case 2.(b).i.A.
     The only frequencies that are separated are $N_4$ and $N_2$ (or, $N_3$ since $N_2 \sim N_3$). So, we use the bilinear refinement of Strichartz on the pair $I_Nu_2 I_Nu_4$ and $L^4$-Strichartz estimate for the other pair.
     \begin{align*}
       Term1 & \lesssim \frac{N_1^3}{m(N_2) m(N_3) m(N_4)} \norm{I_Nu_2
       I_Nu_4}_{L^2(\R^2 \times [0,\delta])} \norm{I_Nu_3}_{L^4(\R^2
       \times [0,\delta])} \norm{I_Nu_1}_{L^4(\R^2 \times [0,\delta])}\\
       & \lesssim \frac{N_1^3}{m(N_2) m(N_3) m(N_4)}
       \frac{N_4^\half}{N_2 \jap{N_2}\jap{N_4}} \frac{1}{\jap{N_3}
       \jap{N_1}} \prod\limits_{i=1}^{4}\norm{I_Nu_i}_{X^\delta_{1,\half+}}\\
       & \lesssim \frac{N_1^2}{m(N_2) m(N_3) m(N_4) N_2^2 N_4^\half
       N_3} \prod\limits_{i=1}^{4}\norm{I_Nu_i}_{X^\delta_{1,\half+}}\\
       & \lesssim \frac{1}{m(N_2) m(N_3) m(N_4) N_2^\half N_3^\half
       N_4^\half} \prod\limits_{i=1}^{4}\norm{I_Nu_i}_{X^\delta_{1,\half+}} 
       \text{ (since } N_1 \lesssim N_2 \sim N_3 \implies N_1^2  \lesssim N_2^\half N_3^\half N_2)\\
       & \lesssim
       N^{-\frac{3}{2}+}N_{max}^{0-}\prod\limits_{i=1}^{4}\norm{I_Nu_i}_{X^\delta_{1,\half+}}.
     \end{align*}
     In the last step, we have used $m(N_i)N_i^p \gtrsim N^p$ for any $p >0$ with $ p+s > 1$, since $N_i \gtrsim N$ for $i=2,3,4$ and $N_{max} \sim N_2$.\\

     \textbf{Case 2.(b).ii.B.} In this case, we have $N \lesssim N_4 \sim N_3 \sim N_2$.\\
     Here, we cannot use the bilinear refinement (\ref{bil}), since no two frequencies  are separated. Therefore, we use $L^4$ Strichartz estimate (\ref{interpolated}) to control the term.  We use the same bound for the multiplier as in Case 2.(b).i.A to get
     \begin{align*}
       Term1 & \lesssim \frac{N_1^3}{m(N_2) m(N_3) m(N_4)}
       \prod\limits_{i=1}^{4} \norm{I_Nu_i}_{L^4(\R^2 \times [0,\delta])}\\
       & \lesssim \frac{N_1^3}{m(N_2) m(N_3) m(N_4)}
       \frac{1}{\jap{N_1}\jap{N_2}\jap{N_3}
       \jap{N_4}}\prod\limits_{i=1}^{4} \norm{I_Nu_i}_{X^\delta_{1,\half+}}\\
       & \lesssim \frac{N_1^2}{m(N_2) m(N_3) m(N_4)} \frac{1}{N_2 N_3
       N_4}\prod\limits_{i=1}^{4} \norm{I_Nu_i}_{X^\delta_{1,\half+}}\\
     & \lesssim \frac{1}{m(N_2) N_2^{\half} m(N_3)N_3^{\frac{1}{4}}
     m(N_4) N_4^{\frac{1}{4}}} \prod\limits_{i=1}^{4} \norm{I_Nu_i}_{X^\delta_{1,\half+}}
     \text{ (since } N_1 \lesssim N_2 \sim N_3 \sim N_4 )\\
     & \lesssim N^{-1+} N_{max}^{0-} \prod\limits_{i=1}^{4} \norm{I_Nu_i}_{X^\delta_{1,\half+}}.
     \end{align*}
     In the last step, we have used $m(N_i)N_i^p \gtrsim N^p$ for any $p > 0$ with $ p + s >1$ since $N_i \gtrsim N$ for $i=2,3,4$ and $N_{max} \sim N_2$.

\begin{rem}
Note that the cases 2.(b).i.B and 2.(b).ii.B provide the worst growth
($N^{-1+}$), which we use in \eqref{stopping} to determine the lower
bound on the Sobolev index, i.e., $\frac{3}{4}$, for which the solution is globally well-posed. Thus, improving the growth in these cases will improve the global well-posedness result.
\end{rem}

Now, we turn to the second term on the right-hand side of \eqref{terms}. 
Again, we perform a dyadic decomposition. Following the previous discussion, we define
$$
Term2 = \abs{\int_{0}^{\delta} \int_{\Gamma_6}
 m(\zeta_1)\, m(\zeta_2)\, m(\zeta_3)\, m( \zeta_4 + \zeta_5 + \zeta_{6})\,(\xi_4+ \xi_5 + \xi_6 + \eta_4 + \eta_5 + \eta_6)
       \prod \limits_{i=1}^{6} \widehat{u_i}(\zeta_i, s)\, ds }.
$$
We arrange the frequencies $N_1, \dots, N_6$ in descending order and
call them $N_1^*,\dots, N_6^*$, respectively (e.g., $N_1^*$ is the
largest frequency, etc.). 
Also, we define $u^*_i = u_i(\zeta_i^*,t)$, where $\zeta_i^*$ is the frequency variable associated with $N_i^{*}$.

For any $1 \le i \le 6$, we have
$$  
|\xi_i + \eta_i| \le |\xi_i|+|\eta_i| \le 2 \sqrt{|\xi_i|^2 + |\eta_i|^2} \lesssim |\zeta_i| \lesssim N_1^*.
$$
Thus,
$$
\abs{ \xi_4 + \xi_5 +\xi_6 + \eta_4 + \eta_5 + \eta_6 } \lesssim N_1^*.
$$
If $ N_1^* << N$, $m(\zeta_1) = m(\zeta_2) = m(\zeta_3) = m(\zeta_4 + \zeta_5 + \zeta_6) = 1$. 
Then, defining $\nu = (1,1)$, we write
$$
Term2  = \abs{\int_{0}^{\delta}
\int_{\Gamma_6}^{}(\zeta_4 + \zeta_5 + \zeta_6) \cdot \nu
\prod\limits_{i=1}^{6} \widehat{u_i}(\zeta_i, s) \, ds}.
$$
Now, permuting the variables $\zeta_1, \dots , \zeta_6$,  we have
\begin{align*}
Term2 & =  \abs{\int_{0}^{\delta} \int_{\Gamma_6}^{}
  \frac{1}{6!} \sum\limits_{\theta \in S_6}^{}
  \left(\zeta_{\theta(4)} + \zeta_{\theta(5)} +
  \zeta_{\theta(6)}\right) \cdot \nu  \prod\limits_{i=1}^{6}
  \widehat{u_i}(\zeta_i, s) ds} \\
  & =  \abs{\int_{0}^{\delta} \int_{\Gamma_6}^{}
  \left(c \sum\limits_{j=1}^{6} \zeta_{j}\right)\cdot \nu
  \prod\limits_{i=1}^{6} \widehat{u_i}(\zeta_i, s) ds} = 0 
\qquad \mbox{for ~some}\quad c \in \N,
\end{align*}
where $S_6$ is the symmetric group on the set $\{1, \dots , 6\}$.

Thus, we assume that $N \lesssim N_1^*$.  
Also, we cannot have $N_1^* >> N_2^* \ge N_3^*$, since the sum of frequencies is zero. Therefore, we have $N_1^* \sim N_2^*$.
We shall break the frequency interactions into the following two cases, based on a comparison between $N_2^*$ and $N_3^*$,
\begin{enumerate}
\item 
$N_1^* \sim N_2^* \sim N_3^* \gtrsim N$
\item 
$N_1^* \sim N_2^* >> N_3^*$ and $ N_1^* \sim N_2^* \gtrsim N$.
\end{enumerate}

\textbf{Case 1.} Using $m(N_1^*)^{2-}(N_1^*)^{2-} \gtrsim N^{2-}$
together with inequalities  \eqref{strsob} and \eqref{alphap}, we get
\begin{align*}
Term2  
  & \lesssim N_1^* \abs{\int_{0}^{\delta} \int_{\Gamma_6}
  \frac{m(\zeta_1^*) m(\zeta_2^*) m(\zeta_3^*) m(\zeta_4^* + \zeta_5^*
  +\zeta_6^*)}{m(\zeta_1^*) m(\zeta_2^*) m(\zeta_3^*) } \prod\limits_{i=1}^{3}
  \widehat{I_Nu_i^*} \prod\limits_{i=4}^{6} \widehat{u_i^*} \, ds}\\
  & \lesssim N_1^* \abs{\int_{0}^{\delta}
  \int_{\Gamma_6}^{}\frac{m(\zeta_2^*)  m(\zeta_3^*)  m(\zeta_4^* + \zeta_5^* +\zeta_6^*)}{m(\zeta_2^*) m(\zeta_3^*)  N_1^*N_2^*N_3^*} \prod\limits_{i=1}^{3}
  \jap{\zeta_i^*}\widehat{I_Nu_i^*} \prod\limits_{i=4}^{6}
  \widehat{u_i^*} \, ds} \\
  & \lesssim \int_{0}^{\delta}
  \int_{\Gamma_6}\frac{1}{(m(N_1^*) N_1^*)^2} \abs{\prod\limits_{i=1}^{3}
  \jap{\zeta_i^*}\widehat{I_Nu_i^*} \prod\limits_{i=4}^{6}
  \widehat{u_i^*}} \, ds 
  \text{ (using } m(\zeta_2^*), m(\zeta_3^*), m(\zeta_4^*+\zeta_5^* +
  \zeta_6^*) \le 1) \\
  & \lesssim \frac{(N_1^*)^{0-}}{N^{2-}} \prod\limits_{i=1}^{3}
  \norm{J^1I_Nu_i^*}_{L^5(\R \times [0, \delta])} \prod\limits_{i=4}^{6}
  \norm{u_i^*}_{L^{\frac{15}{2}}(\R^2 \times [0, \delta])} \\
  & \lesssim \frac{(N_1^*)^{0-}}{N^{2-}} \prod\limits_{i=1}^{3}
  \norm{I_Nu_i^*}_{X^\delta_{1,\half+}}  \prod\limits_{i=4}^{6}
  \norm{u_i^*}_{X^\delta_{\alpha(\frac{15}{2}), \half+}} \\
  & \lesssim \frac{(N_1^*)^{0-}}{N^{2-}} \prod\limits_{i=1}^{6} \norm{I_Nu_i^*}_{X^\delta_{1,\half+}} \\
  & \lesssim N^{-2+} N_{max}^{0-}  \prod\limits_{i=1}^{6} \norm{I_Nu_i^*}_{X^\delta_{1,\half+}} ,
\end{align*}
where we have used (\ref{alphap}) with $\alpha(\frac{15}{2}) =
\frac{1}{3}+ $ and the inequality \eqref{Ismoothing}.
Here, $(\frac{1}{3}+ ) + 1 - s < \frac{7}{12}+ <1$ for any $s \in (3/4,1)$. Thus, we have
\begin{align*}
  \norm{u}_{X^\delta_{\frac{1}{3}+, \half+}} \lesssim \norm{I_Nu}_{X^\delta_{1, \half+}}.
\end{align*}

\textbf{Case 2.}  $N_1^* \sim N_2^* >> N_3^*$ and $N_2^* \gtrsim N$.
We use the following bound
\begin{align*}
\abs{m(\zeta_1)m(\zeta_2) m(\zeta_3) m(\zeta_4 + \zeta_5 + \zeta_6)}
\lesssim 1,
\end{align*}
and estimates \eqref{bil} and \eqref{strsob} to get
\begin{align*}
 & Term2  \lesssim \frac{N_1^*}{m(N_1^*)m(N_2^*)m(N_3^*)}
  \norm{I_Nu_1^* I_Nu_3^*}_{L^2(\R^2 \times [0,\delta])} \norm{I_Nu_2^* \prod\limits_{i=4}^{6} u_i^*}_{L^2(\R^2 \times [0,\delta])} \\
  & \lesssim \frac{N_1^* (N_3^*)^\half}{m(N_1^*)m(N_2^*)m(N_3^*)N_1^*
    \jap{ N_1^* } \jap{N_3^*}}
    \norm{I_Nu_1^*}_{X^\delta_{1,\half+}}\norm{I_Nu_3^*}_{X^\delta_{1,\half+}}
    \norm{I_Nu_2^*}_{L^5(\R^2 \times [0,\delta])}
    \prod\limits_{i=4}^{6} \norm{u_i^*}_{L^{10}(\R^2 \times [0,\delta])} \\
  & \lesssim \frac{1}{m(N_1^*) m(N_2^*)m(N_3^*) N_1^* N_2^* \jap{
  N_3^* }^{\half}}
  \norm{I_Nu_1^*}_{X^\delta_{1,\half+}}\norm{I_Nu_3^*}_{X^\delta_{1,\half+}}
  \norm{I_Nu_2^*}_{X^\delta_{1,\half+}} \prod\limits_{i=4}^{6} \norm{u_i^*}_{L^{10}(\R^2 \times [0,\delta])} \\
  & \lesssim N^{-2+} N_{max}^{0-} \prod\limits_{i=1}^{6}\norm{I_Nu_i^*}_{X^\delta_{1,\half+}}.
\end{align*}
Here, we also used
\begin{align*}
  m(\zeta_1*)(N_1^*)^{1-} \gtrsim N^{1-}, \quad 
  m(\zeta_2*)(N_2^*)^{1-} \gtrsim N^{1-}, \quad 
  m(\zeta_3*)\jap{ N_3^* }^\half \gtrsim 1,
\end{align*}
and 
\begin{align*}
  \norm{u}_{L^{10}(\R^2 \times [0,\delta])} \lesssim
  \norm{u}_{X^\delta_{\alpha(10), \half+}} =
  \norm{u}_{X^\delta_{\half+,\half+}} \lesssim
  \norm{I_Nu}_{X^\delta_{1,\half+}},
\end{align*}
due to \eqref{Ismoothing}, since $(\frac{1}{2}+) + 1 - s < \frac{3}{4}+$  for $ s \in (\frac{3}{4},1)$. This completes the proof of Proposition \ref{propn:1}.
\qed

\section{A variant of the local existence theorem}\label{variant}
Applying the operator $I_N$ on IVP (\ref{gzk}), we get the modified IVP
\begin{equation}\label{igzk}
\begin{cases}
  \partial_t I_Nu + (\partial_x^3 + \partial_y^3) I_Nu + \sigma a (\partial_x+ \partial_y) \left( I_N(u^{3})\right) = 0 \\
  I_Nu(x,y,0) = I_Nu_0(x,y).
\end{cases}
\end{equation}
The Duhamel's formula for the IVP (\ref{igzk}) is
\begin{equation}\label{duhamel}
I_Nu(t) = U(t)I_Nu_0 - \sigma a\int_{0}^{t} U(t-s) (\partial_x+\partial_y) \left( I_N(u^{3}) \right)\, ds.
\end{equation}
To work on $X_{s,b}$ spaces, we consider the following integral equation instead
\begin{equation}\label{iduhamel}
I_Nu(t) = \psi(t)U(t)I_Nu_0 - \sigma a \psi_\delta(t)\int_{0}^{t} U(t-s) (\partial_x+\partial_y) \left( I_N(u^{3}) \right)\, ds,
\end{equation}
where $\psi, \psi_\delta$ are as in Lemma \ref{lemma:duhamel} below.
It is clear that if $u$ is a solution to this equation, then $u \vert_{[0,\delta]}$ is a solution to \eqref{duhamel}.  
We will use the following two Lemmas (from Gr\"unrock and Herr \cite{GH} and Colliander, Keel, Staffilani, Takaoka and Tao \cite{CKSTTKdV})
to establish a local existence result for $I_Nu$ in $H^1(\R^2)$ space.
\begin{lemma}[Lemma 2.2 in \cite{GH}]\label{lemma:duhamel}
Let $\psi \in C_0^{\infty}([-2,2])$ be even, $0\le \psi \le 1 $,  $\psi(t) = 1 $ for $|t|\le 1$. Let $\psi_\delta(t) = \psi(\frac{t}{\delta})$ for $\delta >0$. For all $s,b \in \R$,
\begin{equation}\label{linear}
\norm{\psi U(t) u_0}_{X_{s,b}} \lesssim \norm{u_0}_{H^s}.
\end{equation}
Also, for $-\frac{1}{2} < b' \le 0 \le b \le (b'+1)$ and $0 <\delta <1$,
\begin{equation}\label{nonlinear}
\norm{\psi_\delta \int_{0}^{t} U(t-s) f(s) ds }_{X_{s,b}} \lesssim \delta^{1-b+b'} \norm{f}_{X_{s,b'}}.
\end{equation}
\end{lemma}
\begin{lemma}[Lemma 12.1 in \cite{CKSTTKdV}] \label{lemma:interpolation}
Let $s_0 >0$, $n\ge1$ and $Z, X_1, \dots , X_n$ be translation-invariant Banach spaces. 
If $T$ is a translation-invariant n-linear operator such that
$$
\norm{I_1^s T(u_1,\dots u_n)}_Z \lesssim\prod\limits_{i=1}^n\norm{I_1^s u_i}_{X_i}
$$
for all $u_1, \dots , u_n$ and for all $0 \le s \le s_0$, then,
$$
\norm{I_N^s T(u_1, \dots , u_n)}_Z \lesssim \prod\limits_{i=1}^n \norm{I_N^s u_i}_{X_i}
$$
for all $N \ge 1$, for all $0 \le s \le s_0$, and for all $u_1, \dots , u_n$ with the implicit constant independent of N.
\end{lemma}

\begin{cor}
For all $\frac{3}{4} < s < 1$ and $N >> 1$, we have
\begin{equation}\label{multilinear}
\norm{(\partial_x + \partial_y) I_N(u ^{3})}_{X_{1,-\half++}}
    \lesssim \norm{I_Nu}_{X_{1,\half +}}^{3}.
\end{equation}
\end{cor}
\proof
Using Lemma \ref{lemma:interpolation}, it is enough to show that
\begin{equation}\label{toshowI1}
\norm{(\partial_x + \partial_y) I_1 (u^3)}_{X_{1, -\half++}} \lesssim
\norm{I_1 u}_{X_{1, \half+}}^3.
\end{equation}
Next, note that for any $ b \in \R$ and $\frac{3}{4} < s < 1$ we have
$$  
\norm{I_1f}_{X_{1,b}} \sim \norm{f}_{X_{s,b}}.
$$
Indeed, recalling the definition of $m$ and decomposing the domain of integration in two parts,  we can write
$$  
\norm{I_1f}_{X_{1,b}}^2 = R_1  + R_2,
$$
where
\begin{equation*}
R_1 = \int_{\abs{\zeta} \le 2} \abs{ m_1(\zeta) \jap{\zeta} \jap{\tau - \xi^3 - \eta^3}^b \widetilde{f}(\zeta, \tau) }^2 d\zeta\, d\tau 
\end{equation*}
and
\begin{equation*}
R_2  = \int_{\abs{\zeta} \ge 2} \abs{\frac{\jap{\zeta}}{\abs{\zeta}^{1-s}} \jap{\tau - \xi^3 - \eta^3}^b \widetilde{f}(\zeta, \tau) }^2 d\zeta d\tau.
\end{equation*}
When $\abs{\zeta} \le 2$, we have $m_1(\zeta) \le 1$ and
$$
\jap{\zeta} = \jap{\zeta}^s \jap{\zeta}^{1-s} \lesssim \jap{\zeta}^s,
$$
which implies
$$
R_1 \lesssim \norm{f}_{X_{s, b}}^2.
$$
On the other hand, when $\abs{\zeta} \ge 2$, we have
$$
\jap{\zeta} =  \sqrt{1 + \abs{\zeta}^2} \lesssim \sqrt{2 \abs{\zeta}^2}  \lesssim \abs{\zeta},
$$
and hence,
$$
\frac{\jap{\zeta}}{ \abs{\zeta}^{1-s}} \lesssim \abs{\zeta}^s \le \jap{\zeta}^s.
$$
Thus, we also have
$$
R_2 \lesssim \norm{f}_{X_{s, b}}^2,
$$
and therefore,
$$
\norm{I_1 f}_{X_{1,b}} \lesssim \norm{f}_{X_{s,b}}.
$$
Similarly, we can show that
$$
\norm{f}_{X_{s,b}}  \lesssim \norm{I_1 f}_{X_{1,b}},
$$
since when $\abs{\zeta} \le 1$, we have $\jap{\zeta}^s \lesssim
\jap{\zeta} = \jap{\zeta} m_1(\zeta)$ and when $\abs{\zeta}  \ge 1$,
we have $\jap{\zeta}^s \lesssim \abs{\zeta}^s =
\frac{\abs{\zeta}}{\abs{\zeta}^{1-s}} \lesssim
\frac{\jap{\zeta}}{\abs{\zeta}^{1-s}} \lesssim \jap{\zeta}
m_1(\zeta)$.
Thus, since $I_1$ commutes with $(\partial_x + \partial_y)$, \eqref{toshowI1} is equivalent to \eqref{basic:trilinear}, completing 
the proof.
\qed


Now, we prove a local existence result for the modified IVP \eqref{igzk}. As a consequence, we also obtain a bound on the $X^{\delta}_{1,\half+}$ norm of the modified solution $I_Nu$, uniformly in the existence time, in terms of the initial data.
\begin{thm}\label{thm:xsb}
Let $\frac{3}{4}<s<1$ and $u_0\in H^s(\R^2)$. Then there exists $\delta =
\delta(\norm{I_Nu_0}_{H^1(\R^2)}) >0$ such that the modified IVP \eqref{igzk} has a solution $I_Nu \in C([0,\delta];H^1(\R^2))$ with
\begin{equation}\label{thm:xsb2}
\norm{I_Nu}_{X_{1,\frac{1}{2}+}^\delta} \lesssim \|I_N u_0\|_{H^1}.  
\end{equation}
\end{thm}
\proof
Applying $X_{1,\half+}^\delta$ norm on both sides of \eqref{iduhamel} and applying estimates \eqref{linear} and \eqref{nonlinear} with $ b = \half+$ and $b' = -\half ++$, we obtain
$$
\norm{I_Nu}_{X_{1,\half+}^\delta} \lesssim
\norm{I_Nu_0}_{H^1(\R^2)} + \delta^{\varepsilon} \norm{(\partial_x+\partial_y)(I_N(u^{3}))}_{X_{1,-\half++}^\delta},
$$
where $\varepsilon = (1-b+b') > 0$.
By definition of the localized norm \eqref{def:rest}, we have 
\begin{equation*}
\norm{I_Nu}_{X_{1,\half+}^\delta} \lesssim
\norm{I_Nu_0}_{H^1(\R^2)} + \delta^{\varepsilon} \norm{(\partial_x+\partial_y)(I_N(\phi^{3}))}_{X_{1,-\half++}},
\end{equation*}
where the function $\phi(x,y,t) = u(x,y,t)$ on $\R^2 \times [0,\delta]$ and
\begin{equation}\label{equiv}
\norm{I_Nu}_{X_{1,\half+}^\delta} \sim \norm{I_N\phi}_{X_{1,\half+}}.
\end{equation}
Using the trilinear estimate \eqref{multilinear} for $\phi$ and
the relation \eqref{equiv}, we get
\begin{equation*}
\norm{I_Nu}_{X_{1,\half+}^\delta} \lesssim
\norm{I_Nu_0}_{H^1(\R^2)} + \delta^{\varepsilon}\norm{I_Nu}_{X_{1,\half+}^\delta}^3
\end{equation*}
for some $0 < \delta <1$ and $\varepsilon >0$.
Setting $r=2 \|I_N u_0\|_{H^1}$ and taking $\delta^\epsilon \leq \min\{\frac1{16 \, r^2}, 1 \}$ (and recalling that $\norm{I_N u(t)}_{X_{1,\half+}^\delta}$ is continuous in $\delta$), the operator $I_N$ defined in \eqref{iduhamel} maps the ball $B_r$ in $X_{1,\half+}^{\delta}$ centered at the origin and of radius $r$ into itself and is a contraction by similar argument. Thus, 
$$
\norm{I_Nu}_{X_{1,\half+}^\delta} \lesssim \norm{I_N u_0}_{H^1(\R^2)}, 
$$
completing the proof. 
\qed
\begin{rem}
We can obtain a more precise 
bound on $\delta$ by observing that the trilinear estimate \eqref{basic:trilinear} holds for $b'=0$
(see the proof of Lemma \ref{lemma:trilinear}). Thus, the estimate
\eqref{multilinear} can be modified to
\begin{align*}
  \norm{(\partial_x + \partial_y) I_N (u^3)}_{X_{1,0}} \lesssim
  \norm{I_Nu}_{X_{1, \half+}}^3.
\end{align*}
Moreover, the inhomogeneous linear estimate
\eqref{nonlinear} holds for
$b'=0$. Thus, the proof of the above theorem works for 
$\varepsilon = 1 - b + b' = \half -$, and hence, we get 
\begin{align}\label{delta}
  \delta^{\half -} \sim \frac{1}{\norm{I_Nu_0}_{H^1(\R^2)}^2}.
\end{align}
\end{rem}
%
%

Next we prove an \textit{a priori} bound on the $\dot{H}^1$ norm of the
initial data of the modified IVP \eqref{igzk} in terms of its energy.
In the focusing case, $\sigma = 1$, the size of the $L^2$ norm of the
initial data of \eqref{gzk} has to be bounded from above by the $L^2$
norm of the ground state. 
\begin{lemma}
Let $I_Nu \in C([0,\delta];H^1(\R^2))$ be the solution of the modified IVP \eqref{igzk} given by Theorem \ref{thm:xsb}. 
If $\sigma= - 1$, we have,
\begin{equation}\label{I_Nulessthan1k2}
\norm{\nabla I_Nu_0}^2_{L^2(\R^2)} \lesssim E[I_Nu_0].
\end{equation}
If $\sigma=+1$, the same conclusion (\ref{I_Nulessthan1k2}) holds if 
\begin{equation}\label{smallness}
\norm{u_0}_2 < \sqrt{ 2ab } \, \norm{\varphi}_2. 
\end{equation}
\end{lemma}
\proof 
This lemma is straightforward in the defocusing case, and follows from the Gagliardo-Nirenberg inequality in the focusing case. For convenience of the reader, we include the proof. From \eqref{energy} the energy of $I_Nu$ at time $t=0$ is
\begin{equation*}
E[I_Nu_{0}] = \frac12\int_{\R^2} |\nabla I_Nu_{0}|^2 dxdy - \frac12\int_{\R^2}
(I_Nu_{0})_x(I_Nu_{0})_y \,dx \,dy - \frac{\sigma a}{4} \int_{\R^2} \left( I_Nu_{0} \right)^{4} dx\, dy.
\end{equation*}
When $\sigma= -1$ the last term is positive, thus, 
$$
\norm{\nabla I_Nu_{0}}_{L^2(\R^2)}^2 - \int_{\R^2}^{} (I_Nu_{0})_x (I_Nu_{0})_y dx\,dy \le E[I_Nu_{0}].
$$
The basic inequality $\alpha \beta \leq \frac12(\alpha^2+\beta^2)$ applied to the middle term, yields
$$
\norm{\nabla I_N u_0}_{L^2(\R^2)}^2 \le 
E[I_Nu_{0}] + \frac{1}{2}\norm{\nabla I_Nu_{0}}_{L^2(\R^2)}^2,
$$
or
\begin{equation*}
\norm{\nabla I_Nu_0}^2_{L^2(\R^2)} \lesssim  E[I_Nu_0].
\end{equation*}

Now we turn to the focusing case $\sigma = +1$. It is possible to rewrite the sharp Gagliardo-Nirenberg inequality \eqref{GNcrit} for the symmetrized version with $u$ (the sharp constant will change accordingly), or to go back to the original variable $v$ and use \eqref{GNcrit} as is. Here we use the second approach. Thus,  
using the change of variables \eqref{cov}, we rewrite the multiplier operator as
\begin{align*}
\widehat{I_Nu}(\xi,\eta) 
  & = m(\xi, \eta) \iint_{\R^2} u(x',y') e^{-i(x'\xi + i y' \eta)} dx'\, dy' \\
  & = m(\xi, \eta) \iint_{\R^2} v(x,y)   e^{-i\left( ax(\xi+\eta) + by(\xi - \eta) \right)}\, |2ab|\, dx\, dy \\
  & = |2ab|\, m(\xi, \eta) \widehat{v}(a(\xi+\eta), b(\xi - \eta)).
\end{align*}
Computing the inverse Fourier transform with the change of
variables $p = a(\xi + \eta)$ and $q = b(\xi - \eta)$, we get
\begin{align*}
I_Nu(x',y') 
  & = |2ab| \iint_{\R^2} m(\xi, \eta) \widehat{v}(a(\xi+ \eta) , b(\xi - \eta)) e^{ix' \xi + i y' \eta} d\xi d\eta \\
  & = |2ab| \iint_{\R^2} m\left(\frac{p}{2a} + \frac{q}{2b},
  \frac{p}{2a} - \frac{q}{2b}\right)  \widehat{v}(p,q) e^{i\left( p \frac{x'+y'}{2a} + q \frac{x' - y'}{2b} \right)} \frac{1}{|2ab|} dp dq \\
  & = \widetilde{I}_N v\left(\frac{x'+y'}{2a}, \frac{x' - y' }{2b}\right) 
  = \widetilde{I}_N v (x,y),
\end{align*}
where
$$
\widehat{\widetilde{I}_N f}(\xi, \eta) = \widetilde{m}(\xi,\eta) \widehat{f}(\xi,\eta) \quad \mbox{and} \quad \widetilde{m}(\xi, \eta) = m\left(\frac{\xi}{2a}+ \frac{\eta}{2b}, \frac{\xi}{2a} - \frac{\eta}{2b} \right).
$$ 
In particular, 
\begin{align*}
\widetilde{m}(\xi,\eta) = 
  \begin{cases}
    1 &\text{ if } \frac{\xi^2}{4a^2} + \frac{\eta^2}{4b^2} < N^2\\
    \left(\frac{N}{\frac{\xi^2}{4a^2} + \frac{\eta^2}{4b^2}}\right)^{1-s} &   \text{ if } \frac{\xi^2}{4a^2} + \frac{\eta^2}{4b^2} > 4N^2.
  \end{cases}
\end{align*}
Thus, we have
\begin{equation}\label{Itilde}
I_Nu_0(x',y') = \widetilde{I}_Nv_0(x,y).
\end{equation}
Note that since $I_Nu_0 \in H^1(\R^2)$, so is $\widetilde{I}_Nv_0$.
Applying the sharp Gagliardo-Nirenberg inequality (\ref{GNcrit}) to
$f(x,y) = \widetilde{I}_Nv_0(x,y)$, from (\ref{aenergy}) we get
\begin{align*}
\mathcal{E}[\widetilde{I}_Nv_0]  \ge \half \norm{\nabla
      \widetilde{I}_Nv_0}^2_{L^2(\R^2)}  - \frac{1}{2}
      \frac{\norm{\widetilde{I}_Nv_0}^2_{L^2(\R^2)}}{\norm{\varphi}^2_{L^2(\R^2)}}
      \norm{\nabla \widetilde{I}_Nv_0}^2_{L^2(\R^2)}.
\end{align*}
The condition \eqref{smallness} implies
$$ 
\norm{\widetilde{I}_Nv_0}^2_{L^2(\R^2)} =
\frac{1}{|2ab|}\norm{I_Nu_0}^2_{L^2(\R^2)} \le \frac{1}{|2ab|}
\norm{u_0}^2_{L^2(\R^2)}   < \norm{\varphi}^2_{L^2(\R^2)},
$$ 
and we conclude that
$$ 
\norm{\nabla \widetilde{I}_N v_0}^2_{L^2(\R^2)}  
\lesssim \mathcal{E}[\widetilde{I}_Nv_0].
$$
Next, using the equations \eqref{covNabla} , \eqref{energyDef} and 
\eqref{Itilde}, we can return to $u_0$. Indeed,
$$
\frac{a^2 + b^2 }{|2ab|} \norm{\nabla I_Nu_0}^2_{L^2(\R^2)} -
\frac{2(a^2 - b^2)}{|2ab|} \int_{\R^2}^{} (I_Nu_0)_{x'}
(I_Nu_0)_{y'}(x',y')\, dx'\, dy' \lesssim \frac{a^2 + b^2 }{|2ab|} E[I_Nu_0],
$$
and thus,
$$  
\norm{\nabla I_Nu_0}^2_{L^2(\R^2)} - \frac{2(a^2 - b^2)}{a^2 + b^2} \int_{\R^2}^{} (I_Nu_0)_{x'} (I_Nu_0)_{y'}(x',y')\, dx' \,dy' 
\lesssim E[I_Nu_0].
$$
Recalling that $\dfrac{2|a^2 - b^2|}{a^2 + b^2} = 1$ and splitting the middle term above,  
we conclude that
\begin{equation*}
\norm{\nabla I_N u_0}^2_{L^2(\R^2)} \lesssim E[I_Nu_0],
\end{equation*}
completing the proof. 
\qed

Now, we state a local existence result for the modified rescaled solution
$I_N u_{\lambda}$, which will be used in the proof of the main
theorem \eqref{thm:main}. Under the assumption that the modified
energy of the rescaled solution is uniformly bounded from above, we 
conclude that that the time of existence of $I_N u_{\lambda}$ is a
constant depending only on $\norm{u_0}_{L^2(\R^2)}$ (in particular, independent of the scaling factor $\lambda$).
To arrive at this
conclusion, the assumption \eqref{smallness} is required in the
focusing case.

For $\lambda>0$ define the rescaled solution
\begin{equation*}
u_\lambda(x,y,t)= {\lambda^{-1}} u\left({\lambda^{-1}}x, {\lambda^{-1} y}, {\lambda^{-3} t} \right) \quad \mbox{and} \quad  u_{0,\lambda}(x,y) =  {\lambda^{-1}}u_0( {\lambda^{-1}}x, {\lambda^{-1}}y).
\end{equation*}

\begin{lemma}\label{rema} 
  Assume  that $E[I_N u_{0,\lambda}] < 1$.  Further assume that either
$\sigma = -1$, or $\sigma = +1$ with \eqref{smallness} holds.

Then, there exists $\delta =
\delta\left( \norm{u_0}_2 \right)>0$ such that $I_N u_{\lambda} \in
C\left( [0, \delta]; H^1\left( \R^2 \right) \right)$ with 
\begin{align}
  \norm{I_N u_{\lambda}}_{X_{1,\half+}^\delta} \lesssim 1
  \label{smallerthan1}
\end{align}
If $\sigma=+1$, the same conclusion (\ref{I_Nulessthan1k2}) holds if 
\begin{equation*}
\norm{u_0}_2 < \sqrt{ 2ab } \, \norm{\varphi}_2. 
\end{equation*}
\end{lemma}
\proof
If we assume that $E[I_Nu_{0,\lambda}] < 1$ and
either $\sigma = -1$ or $\sigma = +1$ and \eqref{smallness} holds,
then, in view of \eqref{I_Nulessthan1k2}, we have
\begin{align*}
  \norm{\nabla I_N u_{0,\lambda}}_{L^2(
  R^2)} \lesssim 1.
\end{align*}
Moreover, using $m \le 1$ and the fact that the IVP \eqref{gzk} is
$L^2$-critical,  we have
\begin{align*}
  \norm{I_Nu_{0,\lambda}}_{L^2(\R^2)}\leq
\norm{u_{0,\lambda}}_{L^2(\R^2)} = \norm{u_0}_{L^2(\R^2)}.
\end{align*}
Thus, we have
\begin{align}
\label{specific}
  \norm{I_N u_{0,\lambda}}_{H^1(\R^2)} \lesssim
  \norm{u_0}_{L^2(\R^2)}.
\end{align}
Since $u_\lambda$ solves the IVP \eqref{gzk}, from \eqref{delta}, the time
of existence of $I_N u_{\lambda}$ given by Theorem \ref{thm:xsb} depends only on $\|u_0\|_{L^2}$, that is, $\delta = \delta(\|u_0\|_{L^2})>0$.
Finally, using \eqref{thm:xsb2} and \eqref{specific}, we have \eqref{smallerthan1}.
\qed

\section{Proof of main theorem} \label{mainproof}
As we mentioned in the introduction, it suffices to prove Theorem \ref{thm:main} as it is equivalent to Theorem \ref{thm:mainOriginal}.

\textbf{Proof of Theorem \ref{thm:main}:}
Let $u_0 \in H^s$, where $ \frac{3}{4} < s < 1$. Given any $T >0$, we
will show that the solution $u$ to \eqref{gzk} exists for time
$[0,T]$, which is equivalent to showing that $I_N u_\lambda \in
H^1(\R^2)$ for time $\lambda^3T$. We will do this by iterating Lemma
\ref{rema}.

Note that the $\dot{H}^s$ norm of the rescaled solution is
$\norm{u_{0,\lambda}}_{\dot{H}^s(\R^2)} = \lambda^{-s} \norm{u_0}_{\dot{H}^s(\R^2)}$, and thus, from \eqref{Ismoothing}, we deduce
\begin{equation}\label{H1N}
\norm{I_Nu_{0,\lambda}}_{\dot{H}^1(\R^2)} = N^{1-s}\lambda^{-s} \norm{u_0}_{\dot{H}^s(\R^2)}. 
\end{equation}

From  \eqref{energy}, using a simple bound $-\int_{\R^2}^{}u_xu_y dxdy \le \half \int_{\R^2}^{}|\nabla u|^2dxdy$, Gagliardo-Nirenberg inequality \eqref{GNcrit} and the identity \eqref{H1N}, we get
\begin{align*}
   E[I_Nu_{0,\lambda}] &\le   \norm{\nabla I_Nu_{0,\lambda}}_{L^2(\R^2)}^2 +
   \frac{ a}{4}\norm{I_Nu_{0,\lambda}}_{L^4(\R^2)}^{4} \\
   & \lesssim  \norm{\nabla I_Nu_{0,\lambda}}_{L^2(\R^2)}^2 + \norm{I_Nu_{0,\lambda}}_{L^2(\R^2)}^2 \norm{\nabla I_Nu_{0,\lambda}}_{L^2(\R^2)}^2 \\
   & \lesssim  \norm{\nabla I_Nu_{0,\lambda}}_{L^2(\R^2)}^2\left( 1 + \norm{I_Nu_{0,\lambda}}_{L^2(\R^2)}^2 \right)\\
   & \lesssim N^{2(1-s)} \lambda^{-2s} \norm{u_0}_{\dot{H}^s(\R^2)}^2 \left( 1 + \norm{u_0}^{2}_{L^2(\R^2)} \right)\\
   & \lesssim N^{2(1-s)} \lambda^{-2s}  \left( 1 +
   \norm{u_0}^{2}_{H^s(\R^2)} \right)^2,
 \end{align*}
where we used $\norm{I_Nu_{0,\lambda}}_{L^2(\R^2)}^2 \le 
\norm{u_0}_{L^2(\R^2)}^2$. 
 
Take $\lambda$ such that
\begin{equation*}
C\, N^{2(1-s)}\lambda^{-2s} \left(1 + \norm{u_0}_{H^s(\R^2)}\right)^{2} = \frac{1}{2},   
\end{equation*}
or,
\begin{equation}\label{nlambdarelation2}
\lambda \sim N^{\frac{1-s}{s}}.
\end{equation}
This implies
\begin{align*}
E[I_Nu_{0,\lambda}] \le \frac{1}{2}.
\end{align*}
We apply Lemma \ref{rema} and Proposition \ref{propn:I_Nugrowth} to $u_\lambda$ and conclude that there exists a
$\delta = \delta(\norm{u_0}_{L^2(\R^2)})$  such that 
$$
E[I_N u_\lambda (\delta)]  = E[I_N u_{0,\lambda}] + CN^{-1+}.
$$

Choosing $N$ large,  we have $E[I_N u_\lambda(\delta)] <1$.

Now, since $\norm{I_Nu_{\lambda}(\delta)}_{L^2(\R^2)}\leq \norm{u_{\lambda}(\delta)}_{L^2(\R^2)} = \norm{u_0}_{L^2(\R^2)}$, we can apply Lemma \ref{rema} again with
$t=\delta$ as the starting time, followed by Proposition
\ref{propn:I_Nugrowth}. In other words, starting at $t=\delta$, the
solution $I_Nu_\lambda$ exists for an additional time $\delta =
\delta\left(\norm{u_0}_{L^2(\R^2)}\right)>0$ with
$$
E[I_N u_\lambda (2\delta)]  = E[I_N u_{0,\lambda}] + 2C N^{-1+}.
$$
Note that in the defocusing case (i.e. $\sigma = +1$), the requirement \eqref{smallness} to apply Lemma
\ref{rema} again is trivially satisfied due the mass-conservation law
\eqref{masscons}. Indeed, under the assumption \eqref{smallness}, as
long as the solution exists, we have
\begin{equation}\label{smallnessCheck}
\norm{u(t)}_{L^2(\R^2)} 
= \norm{u_0}_{L^2(\R^2)} < \sqrt{2ab} \norm{\varphi}_{L^2(\R^2)}.
\end{equation}

We repeat this process $M$ times, as long as
$E[I_Nu_\lambda(M\delta)] <1$ and additionally for $\sigma =+1$,  as
long as $\norm{u(M\delta)}_{L^2(\R^2)} <
\sqrt{2ab}\norm{\varphi}_{L^2(\R^2)}$, which holds again by
\eqref{smallnessCheck}. Thus, the iterative application of Lemma
\ref{rema} is valid as long as
\begin{equation}\label{stopping}
MCN^{-1+} < \frac{1}{2} \implies M \lesssim N^{1-}.
\end{equation}

 To show that the solution $I_Nu_\lambda$ exists for time
 $\lambda^3T$, we need that
$$
M\delta > \lambda^3T \implies N^{1-} \gtrsim \lambda^3T,
$$
where we have used (\ref{stopping}) and the fact that $\delta$ depends only on $\norm{u_0}_{L^2(\R^2)}$.
Using the relation (\ref{nlambdarelation2}) between $N$ and $\lambda$ , we have
\begin{align*}
N^{1-} \gtrsim \lambda^3T \sim N^{\frac{3(1-s)}{s}} T
\implies N^{\left[  1-\frac{3(1-s)}{s}\right]-} > cT .
\end{align*}
Now, we need the  power of $N$ to be positive so that $T$ can be taken as large as we want. Therefore,
$$
1 - \frac{3(1-s)}{s}  >  0 \iff  s > \frac{3}{4}.
$$
Finally, we derive a polynomial bound for the $H^s$ norm of the solution to IVP \eqref{gzk}. Note that in the previous argument we can select time $ T \sim N^{\frac{4s - 3}{s}-}$.
By definition,
\begin{equation}\label{hsnorm}
\norm{u(T)}_{H^s(\R^2)}^2  \lesssim \norm{I_Nu(T)}_{H^1(\R^2)}^2  =
\norm{I_Nu(T)}_{L^2(\R^2)}^2 + \norm{I_N u(T)}_{\dot{H}^1(\R^2)}^2 .
\end{equation}
The first term on the right hand side can be bounded by $\norm{u_0}_{L^2}^2$, since $m \le 1$ and mass conservation \eqref{masscons} holds.  
To bound the second term on the right-hand side of \eqref{hsnorm}, first, we note that $E(I_Nu_\lambda)(\lambda^3T) \lesssim 1$.
Moreover, when $\sigma =+1$, for any time $t \in [0,T]$, we have $\norm{u_\lambda(t)}_{L^2(\R^2)} \lesssim \sqrt{2ab} \norm{\varphi}_{L^2(\R^2)}$ from \eqref{smallnessCheck}.
Then using relations \eqref{I_Nulessthan1k2} and \eqref{nlambdarelation2},
we get
\begin{align*}
\norm{I_N u(T)}_{\dot{H}^1(\R^2)}^2 & = \lambda^2 \norm{I_N
u_\lambda(\lambda^3T)}_{\dot{H}^1(\R^2)}^2
   \lesssim  \lambda^2 E[I_Nu_\lambda](\lambda^3T)
   \lesssim  \lambda^2
   \lesssim N^{\frac{2(1-s)}{s}} 
   \sim  T^{\frac{2(1-s)}{4s-3}+}.
\end{align*}
From (\ref{hsnorm}), we have
$$
\norm{u(T)}_{H^s(\R^2)}^2 \lesssim \left( 1 + T \right)^{\frac{2(1-s)}{4s-3}+},
$$
which concludes the proof of Theorem \ref{thm:main}.
\qed


\end{document}